\documentclass[10pt,a4paper]{amsart}
\usepackage{amsmath,amsfonts,amsthm}
\newtheorem{thm}{Theorem} [section]
\newtheorem{cor}[thm]{Corollary}
\newtheorem{lem}[thm]{Lemma}
\newtheorem{prop}[thm]{Proposition}
\theoremstyle{definition}
\newtheorem{defn}[thm]{Definition}
\theoremstyle{theorem}
\newtheorem{rem}[thm]{Remark}
\newtheorem{fact}[thm]{Fact}

\newcommand{\G}{\Gamma}
\renewcommand{\>}{\rangle}
\newcommand{\<}{\langle}
\usepackage{color}

\date{}
\begin{document}

\author{A. Ould Houcine}
\title{Homogeneity and prime models in  torsion-free hyperbolic groups}
\maketitle

\begin{abstract}  We show that any nonabelian free group $F$ of finite rank is homogeneous; that is for any  tuples $\bar a$, $\bar b \in F^n$, having the same complete $n$-type,  there exists an automorphism of $F$  which sends $\bar a$ to $\bar b$. 

We further study existential types and we show that for any  tuples $\bar a, \bar b \in F^n$,  if $\bar a$ and $\bar b$ have the same  existential $n$-type, then either $\bar a$ has  the same existential type as  a power of a primitive element, or  there exists an existentially closed subgroup $E(\bar a)$ (resp. $E(\bar b)$) of $F$ containing $\bar a$ (resp. $\bar b$) and an isomorphism    $\sigma : E(\bar a)  \rightarrow E(\bar b)$   with $\sigma(\bar a)=\bar b$.

We will deal with  non-free two-generated torsion-free hyperbolic groups and we show that they are $\exists$-homogeneous and prime.  This gives,   in particular,  concrete examples of finitely generated groups which are prime and not QFA. 
\end{abstract}

\section{Introduction}

From a model-theoretical point of view, the homogeneity can be seen as a kind of saturation.  For instance,  a countable model in a countable language is saturated,  if and only if, it is homogeneous and realizes all types of its complete theory. Homogeneity  is also a notion related to prime  models and it is well known that a countable prime model in a countable language  is homogeneous.

It is easy to see that a free group is not saturated. Consequently, it is natural to wonder if any free group is at least homogeneous.  This question was studied in the case of the free group of rank $2$ in \cite{nies-free}, where A. Nies proved that this last group is $\exists$-homogeneous and not prime. 

In this paper, we study the homogeneity of free  groups  of higher rank and that of particular torsion-free hyperbolic groups including the  two-generated ones.  The study of these last groups was widely motivated by the previous result of A. Nies, where the proof seems to use strongly the two-generation property. We emphasize  that, 
by a result of T. Delzant [Del96], any (torsion-free) hyperbolic group is embeddable in a two-generated (torsion-free) hyperbolic group. In some  sense, these last groups can have a very complicated structure.

Let $\mathcal M$ be a  model, $P$ a subset of $\mathcal M$ and $\bar a$ a tuple from   $\mathcal M$. The  \textit{type} (resp. \textit{existential type}) of $\bar a$ over $P$, denoted $tp^{\mathcal M}(\bar a|P)$ (resp. $tp_{\exists}^{\mathcal M}(\bar a|P)$), is the set of formulas $\varphi (\bar x)$ (resp. existential formulas $\varphi(\bar x)$) with parameters from $P$ such that $\mathcal M$ satisfies $\varphi(\bar a)$.   

A countable  model $\mathcal M$  is called  \textit{homogeneous} (reps. \textit{$\exists$-homogeneous}), if for   any tuples $\bar a$, $\bar b$ of $\mathcal M^n$,  if $tp^{\mathcal M}(\bar a)=tp^{\mathcal M}(\bar b)$ (resp. $tp_{\exists}^{\mathcal M}(\bar a)=tp_{\exists}^{\mathcal M}(\bar b)$) then there exists an automorphism  of $\mathcal M$ which sends $\bar a$ to $\bar b$.  We note,  in particuliar, that $\exists$-homogeneity implies homogeneity.  For further  notions of homogeneity,  we refer the reader to \cite{Hodges(book)93, Marker}.  

We recall also that a model $\mathcal M$ is a said to be  \textit{prime},  if it is elementary embeddable in every model  of  its  complete theory.  As usual, to axiomatize  group theory,  we use the language $\mathcal L=\{., ^{-1}, 1\}$, where $.$ is interpreted by the multiplication, $^{-1}$ is interpreted by the function which sends every element to its inverse and $1$ is interpreted by the trivial element.  The main results of this paper are as follows.  

\begin{thm} \label{thm1} Let $F$ be a nonabelian free group of finite rank. For any  tuples $\bar a, \bar b \in F^n$ and for any subset $P \subseteq F$,  if $tp^F(\bar a|P)=tp^F(\bar b|P)$ then there exists an automorphism   of $F$  fixing pointwise $P$ and sending $\bar a$ to $\bar b$.
\end{thm}

Let $\mathcal M$ be a model and $\mathcal N$ a submodel of $\mathcal M$. The model $\mathcal N$ is said to be \textit{existentially closed} (abbreviated e.c.)  in $\mathcal M$,  if for any existential formula $\varphi(\bar x)$ with parameters from $\mathcal N$, if $\mathcal M \models \exists \bar x \varphi(\bar x)$,  then $\mathcal N \models \exists \bar x \varphi(\bar x)$.

\begin{defn}  Let $F$ be a free group and let $\bar a=(a_1, \dots, a_m)$ be a tuple from $F$.  We say that $\bar a$ is \textit{a power of a primitive element} if there exist integers $p_1, \dots, p_m$ and a primitive element $u$ such that $a_i=u^{p_i}$ for all $i$. 
\end{defn}

\begin{thm}\label{thm2} Let $F$ be a nonabelian free group of finite rank. Let  $\bar a, \bar b \in F^n$ and  $P \subseteq F$ such that $tp_{\exists}^F(\bar a|P)= tp_{\exists}^F(\bar b|P)$. Then either $\bar a$ has the same  existential type as a power of a primitive element,  or    there exists an existentially closed subgroup $E(\bar a)$ (resp. $E(\bar b)$) containing $P$ and $\bar a$ (resp. $\bar b$)   and an isomorphism $\sigma : E(\bar a)  \rightarrow E(\bar b)$  fixing pointwise $P$ and sending $\bar a$ to $\bar b$.   
\end{thm}

A group $\G$ is said \textit{co-hopfian}, if any injective endomorphism of $\G$ is an automorphism. In \cite{Sela-hopf}, Z. Sela proved  that a non-cyclic freely indecomposable torsion-free hyperbolic group is co-hopfian.  When the given group  is two-generated, we have  in fact a more strong property.  We introduce the following definition.

\begin{defn}  \label{def-strongly-co} A group  $\G$ is said to be \textit{strongly co-hopfian}, if there exists a finite subset $S \subseteq \G\setminus \{1\}$  such that for any endomorphism $\varphi$ of $\G$, if $1 \not \in  \varphi(S)$ then $\varphi$ is an automorphism.  \qed
\end{defn}

For instance,   Tarski monster groups   are strongly co-hopfian.  Recall that a Tarski monster  group is an infinite group in which every nontrivial proper subgroup is cyclic  of order $p$,  where $p$ is a fixed prime.   Such groups were built by A. Ol'shanski{\u\i} in \cite{Ol-monster} and for more details we refer the reader to \cite{book-olshanski}. It is easily seen that they are simple.  It is an immediate consequence that a nontrivial  endomorphism of a Tarski Monster group is an automorphism and thus a such group is  strongly co-hopfian.

\begin{thm} \label{thm3} A non-free two-generated torsion-free hyperbolic group is strongly co-hopfian. 
\end{thm}

The proof of Theorem \ref{thm3} is related to   properties of $\G$-limit groups and   to special properties of two-generated hyperbolic groups. We will also use   the following   notion.

\begin{defn} \cite[Definition 3.4]{Ould-equa}
A finitely generated   $\G$-limit group $G$ is said
\emph{$\G$-determined} if  there exists a finite subset $S
\subseteq G \setminus \{1\}$ such that for any homomorphism $f : G
\rightarrow L$, where $L$ is a $\G$-limit group, if $1 \not \in
f(S)$ then $f$ is an embedding.  \qed
\end{defn} 

From Theorem \ref{thm3},  we  deduce  the following. 

\begin{cor} \label{cor1}   A non-free   two-generated torison-free hyperbolic group $\G$ is   $\exists$-homogeneous, prime and $\G$-determined. 
\end{cor}

The above enables one to give examples of one-relator $\exists$-homogeneous and prime groups. Indeed, in the free group $F=\<a,b|\>$  if we let  $r \in F$ such that $r$ is  root-free  and the symmetrized set that it generates satisfies the small cancellation condition $C'(1/6)$,  then  the group $\G=\<a,b| r=1\>$ is a non-free two-generated  torsion-free one-relator hyperbolic group, which is consequentely $\exists$-homogeneous and prime.

\textit{Rigid groups} are defined in \cite{rips-sela-structure-rigidity} and an equivalent definition in our context is that  a torsion-free hyperbolic  group $\G$   is  \textit{rigid} if it is freely indecomposable and does not admit an essential cyclic splitting.    Following  \cite[Definition 7.1]{groves-2007},   a finitely generated subgroup $H\leq \Gamma$ is called \textit{immutable} if there are finitely many embeddings $\varphi_1, \dots, \varphi_n: H \rightarrow \Gamma$ so that any embedding  $\varphi : H \rightarrow \Gamma$   is conjugate to one of the $\varphi_i$. It follows by \cite[Lemma 7.2]{groves-2007} that a subgroup $H\leq \G$ is immutable if and only if it is rigid. We note that a torsion-free hyperbolic group is an immutable subgroup of  itself  if and only if it is rigid.  We point out  that a rigid torsion-free hyperbolic group is strongly co-hopfian and in particular $\exists$-homogeneous and prime, as well as co-hopfian immutable subgroups of torsion-free hyperbolic groups (see Lemma \ref{main-prop}, Lemma \ref{main-lem2} and \ref{main-lem3}). 

As it was mentioned  in \cite{groves-2007} (see \cite{Beleg}), the fundamental group  of a closed hyperbolic $n$-manifold where $n \geq  3$ is rigid and thus it is $\exists$-homogeneous and prime. Hence this gives examples of $\exists$-homogeneous torsion-free hyperbolic groups which are not necessarly two-generated.

We notice that Corollary  \ref{cor1} shows also that the Cantor-Bendixson rank of  a two-generated  torsion-free hyperbolic group is $0$ in the space of its limit groups. For more details on this notion, we refer the reader to \cite{ould-cb-limit, ould-rank}.

Following \cite{Nies-QFA},  a finitely generated group $\G$ is said to be QFA (for quasi-finitely axiomatizable) if there exists a sentence $\varphi$ satisfied by $\G$ such that any finitely generated group satisfying $\varphi$ is isomorphic to $\G$.  A. Nies \cite{Nies-QFA} has proved that the free nilpotent group of class 2 with two generators is QFA prime. F. Oger and G. Sabbagh show that finitely generated nonabelian free nilpotent groups are QFA and prime \cite{Oger-Sabbagh}. It is proved in  \cite{Nies-QFA2}  the existence of  countinousely  many non-isomorphic finitely generated prime groups,  which  implies that there exists a finitely generated  group which is a prime  but not  QFA.

Corollary \ref{cor1} gives  concrete examples of finitely generated groups which are prime  and which are not QFA.  Indeed, it follows from \cite{sela-hyp}, that if $\G$ is a  non-free two-generated torsion-free hyperbolic group,  then $\G$ is an elementary subgroup of $\G * \mathbb Z$; and thus $\G$ is not QFA.

The present paper is organized as follows. In the next section, we record the material that 
we require around $\G$-limit groups  and  tools needed in the sequel. In Section 3,   we show preliminary propositions.  Section 4 concerns  existential types and the proof of Theorem \ref{thm2}  when $P=\emptyset$.  Section 5 is devoted to the general case and we show Theorem \ref{thm1}  when $P=\emptyset$. Section 6  deals with parameters and the proof of Theorem \ref{thm1} and \ref{thm2} for any $P$.   Section 7 is devoted to the proof of  Theorem \ref{thm3} and Corollary \ref{cor1}.  Section 9 
concludes with some remarks and we show, in particular, that non-cyclic torsion-free hyperbolic groups are connected.

\begin{rem} When the work presented in this paper was under verification and more thorough investigation, the preprint \cite{perin-hom} appeared where C. Perin and  R. Sklinos show  the homogeneity of  countable free groups and give a counter-example in the case of  torsion-free  hyperbolic groups. The method which we use in this paper is  different from that used in \cite{perin-hom}.  \end{rem}

\noindent\textbf{Acknowledgements}.  The author wishes to express his thanks to G. Sabbagh and A. Ivanov  for suggesting several questions.

\section{Prerequisites}

We recall some material  about  $\G$-limit groups, where $\G$ is a torison-free hyperbolic group, developped by Z. Sela \cite{sela-hyp}. For more details, we refer the reader to  \cite{sela-hyp}; see also \cite{groves-2007, chloe-these}.  We begin by giving the basic definitions. 

\begin{defn} Let $\G$ be a group. 

\smallskip
$(1)$ A sequence of homomorphisms $(f_n)_{n \in \mathbb N}$ from  $H$ to $\Gamma$ is called \textit{stable} if for any $h \in H$ either $f_n(h)=1$ for all but finitely many  $n$, or   $f_n(h) \neq 1$ for all but finitely many $n$.  The \textit{stable 
kernel} of $(f_n)_{n \in \mathbb N}$, denoted $Ker_{\infty}(f_n)$,  is the set of elements $h \in H$ such that $f_n(h)=1$ for all but finitely many $n$. 

\smallskip
$(2)$  A \textit{$\Gamma$-limit group} is a group $G$ such that  there exists a group  $H$  and a stable sequence of homomorphisms $(f_n)_{n \in \mathbb N}$ from $H$ to $\G$ such that $G=H/ Ker_{\infty}(f_n)$.
\end{defn}

Let $G$ be a group and $A$ a subgroup of $G$.  The group $G$ is said \textit{freely $A$-decomposable} or \textit{freely decomposable relative to A},  if $G$ has a nontrivial free decomposition $G=G_1*G_2$ such that $A \leq G_1$.  Otherwise, $G$ is said freely $A$-indecomposable or freely indecomposable relative to $A$. 

A \textit{cyclic splitting} of a group $G$ is a graph of groups decomposition of $G$ in which every edge group is infinite cyclic.  A cyclic splitting is said \textit{essential} if any edge group has infinite index in the adjacent vertex groups.

\begin{thm} \emph{\cite{sela-hyp}(see also \cite[Theorem 3.9]{groves-2007} \label{thm-princip})}Let $\G$ be a torsion-free hyperbolic group. Let $H$ be a freely indecomposable finitely generated group and let $(f_n : H  \rightarrow \G)_{n \in \mathbb N}$ be  a stable sequence of pairwise nonconjugate  homomorphisms with trivial stable kernel. Then $H$ admits an essential cyclic splitting. \qed
\end{thm}

Z. Sela \cite{Sela-Dio-VI} and O. Kharlampovich and A. Myasnikov \cite{Kharlam-Mya-free-gps} show that nonabelian free groups have the same elementary theory. More exactly, the following stronger result. 

\begin{thm} \label{thm-elemantary-free-factor} A nonabelian free factor of a free group of finite rank is an elementary subgroup. \qed
\end{thm}

In \cite{Sela-Dio-V1, Sela-Dio-V2}, Z. Sela shows the following quantifier-elimination result. 

\begin{thm} Let $\varphi(\bar x)$ be a formula. Then there exists a boolean combination  of $\exists \forall$-forumula $\phi(\bar x)$,  such that for any nonabelian free group $F$ of finite rank, one has $F \models \forall \bar x(\varphi(\bar x) \Leftrightarrow \phi(\bar x))$. \qed
\end{thm}

We notice, in particular, that if $\bar a, \bar b \in F^n$ such that $tp_{\exists \forall}^F(\bar a)=tp_{\exists \forall}^F(\bar b)$,  then $tp^F(\bar a)=tp^F(\bar b)$. 

In \cite{chloe-these}, the converse of Theorem \ref{thm-elemantary-free-factor} is proved. 

\begin{thm} \label{thm-free-factor}An elementary subgroup of a free group of finite rank is a free factor. \qed
\end{thm}

In \cite{pillay-primitive}, A. Pillay shows the following. 

\begin{thm} \label{primitive} Let $F$ be a nonabelian free group of finite rank and $u,v \in F$ such that $tp^F(u)=tp^F(v)$. If $u$ is primitive, then $v$ is primitive. \qed
\end{thm}

In the proof of Theorem \ref{thm1} and  \ref{thm2},  we use the following properties of free groups.  Let $F$ be a free group with a finite basis $A$. Let  $|u|$ denote the length of a word $u$ in $F$, with respect to the basis $A$.  From \cite[Proposition 2.5, Ch I]{LyndonSchupp77}, a subgroup $H \leq F$ has a Nielsen-reduced generating set $U$ and a Nielsen-reduced set $U$ satisfies the following property \cite[Proposition 2.13, Ch I]{LyndonSchupp77}:  if $w \in H$ has the form $w=u_1 u_2 \dots u_m$ where each $u_i\in U^{\pm 1}$ and $u_i u_{i+1} \neq 1$ then $|w|\geq m$ and $|w|\geq |u_i|$ for any $i$.
Hence we can conclude the previous remarks with the following proposition.

\begin{prop}  \label{cor2} Let $F$ be a nonabelian free group. Then any subgroup $H$   of rank $m$ of $F$ has a basis $B=\{b_1, \dots, b_m\}$ such that for any reduced nontrivial word $w$ on $A$  one has $|w| \geq |b|$ for any $b \in B$ which appears in the reduced form of $w$ with respect to $B$.  \qed 
\end{prop}

\begin{prop} \label{prop-hom-gp-libre} \cite[Proposition 2.12]{LyndonSchupp77} Let $f$ be  homomorphism from a free group $F$ of finite rank onto a free group $G$. Then $F$ admits a free decomposition $F=A*B$ such that $f(A)=G$ and $f(B)=1$ and $f$ is injective on $A$. \qed
\end{prop}

The next proposition is a particular case of \cite[Proposition 1]{john-turner}.

\begin{prop} \label{prop1} Let $F$ be a nonabelian free group of finite rank and let $H$ be a  nontrivial subgroup of $F$. If $f : F \rightarrow F$ is a non-surjective monomorphism such that $f(H)=H$, then $F$ is freely $H$-decomposable. \qed
\end{prop}

In dealing with prime models, the following characterization is useful. 

\begin{prop}\label{prop-prime-model} \cite{Hodges(book)93} Let $\mathcal M$  be a countable model. Then $\mathcal M$ is a prime model 
of its theory if and only if  for every $m \in \mathbb N$, each orbit under the action of $Aut(\mathcal M)$ on $\mathcal M^m$ 
is first-order definable without parameters. 
\end{prop}

We end this section with the following theorem needed elsewhere. 
\begin{thm} \emph{\cite[Theorem 1.22]{sela-hyp}} \label{thm-sela} Any system of equations  in finitely many variables is equivalent in a trosion-free hyperbolic group to a finite subsystem. \qed
\end{thm}

It follows that a torsion-free hyperbolic group is equationnally noetherian. For more details on this notion, we refer the reader to \cite{Baum-Rem-algebraic}.

\section{Preliminaries}

Recall that a subgroup $A$ of a group $G$ is said to be \textit{malnormal} if for any $g \in G\setminus A$, $A^g \cap A=1$. A group $G$ is said a \textit{CSA-group}, if any maximal abelian subgroup of $G$ is malnormal.  In particular a CSA-group is commutative transitive; that is the commutation is a transitive relation on the set of nontrivial elements. Basic facts  about CSA-groups and their HNN-extensions will be used freely through the rest of the paper.  For more details, see \cite{Ould-csa-super, ould-e.c.groups, ould-csa-subgroup}.  In an HNN-extension we denote by $|g|$ the length of normal forms of $g$.

\begin{lem} \label{lem1}  Let $G=\<H,t|U^t=V\>$ where $U$ and $V$ are cyclic subgroups of $G$ generated respectively by $u$ and $v$. Suppose that:

$(i)$  $U$ and $V$ are   malnormal in $H$. 

$(ii)$ $U^h \cap V =1$ for any $h \in H$.

\smallskip
Let $\alpha, \beta \in H$, $s \in G$  such that $\alpha ^s=\beta$,  $|s| \geq 1$. Then one of the following cases holds: 

\smallskip
$(1)$   $\alpha= {u^p}^{\gamma}$,     $\beta= {v^p}^{\delta}$,  $s=\gamma^{-1} t \delta$,  where $p \in \mathbb Z$ and $\gamma, \delta \in H$. 

\smallskip
$(2)$   $\alpha= {v^p}^{\gamma}$,     $\beta= {u^p}^{\delta}$,  $s=\gamma^{-1} t^{-1} \delta$,  where $p \in \mathbb Z$ and $\gamma, \delta \in H$. 
\end{lem}

\proof   Write $s=h_0t^{\epsilon_0}\cdots t^{\epsilon_n}h_{n+1}$ in normal form.  Hence 
$$
h_{n+1}^{-1}t^{-\epsilon_n}\cdots t^{-\epsilon_0}h_0^{-1} \alpha h_0t^{\epsilon_0}\cdots t^{\epsilon_n}h_{n+1}=\beta,
$$
and thus either $h_0^{-1} \alpha h_0 \in U$  and $\epsilon_0=1$ or $h_0^{-1} \alpha h_0 \in V$ and $\epsilon_0=-1$.

We treat only the first case, the other case can be treated similarly. Therefore   $\alpha = h_{0}u^ph_0^{-1}$  for some $p \in \mathbb Z$.

We claim that $n=0$. Suppose that $n \geq 1$. Then $h_1^{-1} v^p h_1 \in U$  and $\epsilon_1=1$ or $h_1^{-1} v^p h_1 \in V $ and $\epsilon_1=-1$.  Since $U^h \cap V=1$, the first case is impossible. Therefore we get the second case and thus  $h_1 \in V$ by the malnormality of $V$.  Hence the sequence $(t^{\epsilon_0}, h_1, t^{\epsilon_1})$ is not reduced; a contradiction.  Thus $n=0$ as claimed and hence  
$
\alpha = h_{0}u^ph_0^{-1}, s= h_0th_1, \beta= h_1^{-1}v^ph_1
$. \qed

\smallskip

In \cite{kapovich-two} the structure of two-generated torsion-free hyperbolic groups was investigated. The next theorem is a version  of \cite[Proposition 5.3]{kapovich-two}. The proof is essentially the same. We remove the occurence of free products with amalgamation and we show  that the cyclic  subgroups involved can be shosen  malnormal in the vertex group. 

\begin{thm} \label{thm-structure-two}Let $\G$ be a non-free two-generated torsion-free hyperbolic group.  Then there exists a sequence of subroups $\G=\G_1 \geq \G_2 \geq \dots \geq \G_n$ satisfying the following properties: 

\smallskip
$(i)$ Each $\G_i$ is two-generated, hyperbolic and quasiconvex; 

$(ii)$ $\G_{i}=\<\G_{i+1}, t | A^t=B\>$, where $A$ and $B$ are a nontrivial malnormal cyclic subgroups of $\G_{i+1}$;

$(iii)$ $\G_n$ is a rigid subgroup of $\G$. 
\end{thm}

\proof 

The construction of the sequence proceeds as follows. If $\G$ is rigid, then we get our sequence by setting $n=1$. Se we  suppose that $\G$ is not rigid. Hence $G$ admits an essential cyclic splitting. By \cite[Proposition 4.26]{Champ-Guirardel}, $\G$ admits an abelian splitting in which each edge group is maximal in the neighbouring vertex group. Since  $\G$ is a trosion-free byperbolic group, any abelian subgroup is cyclic and thus  $\G$ admits a cyclic splitting in which each edge group is malnormal in the neighbouring vertex group. Since that splitting is cyclic  and essential  it  has only one vertex and one edge by \cite[Theorem A]{kapovich-two}. Hence, $\G$ admits a splitting of the form  
$$
\G=\<K, t |A^t=B\>,
$$
where $A$ and $B$ are a malnormal nontrivial  cyclic subgroups of $K$.  By \cite[Proposition 3.8]{kapovich-two}, $ K$ is two-generated and freely indecomposable. It is also  quasiconvex by \cite[Proposition 4.5]{kapovich-two} and hence hyperbolic. We set $\G_{2}=K$. 

Now $\G_2$ satisfies the same properties as $\G$ and we can apply the same procedure to it as above. As in the proof of \cite[Proposition 5.3]{kapovich-two}, we show that there exists $p \in \mathbb N$ so that any sequence of subroups $\G_1=\G \geq \G_2 \geq \dots \geq \G_n$ satisfying $(i)$-$(ii)$ of the theorem we have $n \leq p$.  By a result of T. Delzant \cite{delzant-two}, the number of distinct conjugacy classes of two-generated freely indecomposable subgroups in a torsion-free hyperbolic group is finite.  Let $p$ to be  that number.  Suppose that $n>p$. Then there exists $i<j<n$ such that $\G_j=\G_i^g$ for some $g \in \G$. Therefore $\G_i^g$ is a proper subgroup of $\G_i$; a contradiction with \cite[Lemma 4.6]{kapovich-two} as $\G_i$ is quasiconvex subgroup of $\G$.

Hence in a maximal sequence $\G=\G_1 \geq \G_2 \geq \dots \geq \G_n$ satisfying $(i)$-$(ii)$ of the theorem, $\G_n$ does not admit an essential splitting and thus it is rigid.\qed

\smallskip
Since  rigid  trosion-free hyperbolic groups are freely indecompoable, their $\exists$-homogeneity and primeness   is a consequence of  the following  lemmas of independent interest.  

\begin{lem}\label{main-prop}
Let $\G$ be a torsion-free hyperbolic group. A rigid   finitely generated subgroup of $\G$ is $\G$-determined. 
\end{lem}

\proof Let $H$ be a finitely generated rigid subgroup of $\G$ and suppose for a contradiction that $H$ is not $\G$-determined. Therefore for any finite subset $S \subseteq H\setminus \{1\}$,   there exists a non-injective  homomorphism $\varphi : H \rightarrow L$, where $L$ is a $\G$-limit group, such that $1 \not \in \varphi(S)$. Since $L$ is a $\G$-limit group, we may suppose without loss of generality that $\varphi : H \rightarrow \G$. Write $H\setminus \{1\}$ as  an increasing sequence of finite subsets $(S_i)_{i \in \mathbb N}$.    Thus there exists a  sequence of non-injective  homomorphisms $\varphi_i : H \rightarrow \G$ such that $1 \not \in \varphi_i(S_i)$. Clearly such a sequence   is stable and   has a  trivial stable  kernel. Since each $\varphi_i$ is non-injective, we can extract a stable subsequence of pairwise nonconjugate homomorphisms with trivial stable kernel.   Hence $H$ admits  an essential  cyclic splitting by Theorem \ref{thm-princip}, which is a contradiction. \qed

\smallskip
We note that a co-hopfian $\G$-determined $\G$-limit group is strongly co-hopfian.

\begin{lem}\label{main-lem2}
Let $\G$ be a torsion-free hyperbolic group. A co-hopfian $\G$-determined $\G$-limit group is $\exists$-homogeneous. 
\end{lem}

\proof 

 Let $H$ be a co-hopfian finitely generated  $\G$-determined $\G$-limit group. Let 
$$
H=\<\bar x | w_i(\bar x )=1, i \in \mathbb N\>
$$

be a presentation of $H$. Since $H$ is $\G$-detremined, there exists a finite number of words $v_1(\bar x )$, $\dots$, $v_m(\bar x)$ such that 
$$
H \models v_i(\bar x) \neq 1,
$$ and for any $\G$-limit group $L$,  if $f : H \rightarrow L$ is a homomorphism such that $v_i(f(\bar x)) \neq 1$ for every $1 \leq i \leq m$, then $f$ is an embedding. 

By Theorem \ref{thm-sela},   there exists $p \in \mathbb N $ such that 
$$
\G \models \forall \bar x ( w_1(\bar x)=1 \wedge \dots \wedge w_p(\bar x)=1 \Rightarrow w_i(\bar x)=1), \leqno (1)
$$
for any $i \in \mathbb N$. Since $H$ is embeddable in  $\G$,  $H$ satisfies the  sentence appearing in $(1)$. 

We conclude that for any tuple  $\bar y$ in $H$ which satisfies 
$$
w_1(\bar y)=1 \wedge \dots \wedge w_p(\bar y)=1 \wedge v_1(\bar y)\neq 1 \wedge \dots  \wedge v_m(\bar y)\neq 1, 
$$
there is a homomorphism $f : H \rightarrow H$ which sends $\bar x$  to $\bar y$ and such a homomorphism is necessarly a monomorphism. Furtheremore, it is an automorphism as $H$ is co-hopfian.

Let $\bar a, \bar b$ be tuples of $H$ such that $tp_{\exists}(\bar a)\subseteq tp_{\exists}(\bar b)$ and let us show that there is an automorphism of $H$ which sends $\bar a$ to $\bar b$. 

Let $\bar a=u(\bar x)$. Since $tp_{\exists}(\bar a)\subseteq tp_{\exists}(\bar b)$, there exists $\bar y$ such that 
$$
\bar b=u(\bar y), w_1(\bar y)=1 \wedge \dots \wedge w_p(\bar y)=1 \wedge v_1(\bar y)\neq 1 \wedge \dots  \wedge v_m(\bar y)\neq 1,
$$
and thus there is an automorphism of $H$ which sends $\bar a$ to $\bar b$. \qed

\smallskip
We introduce the following definition,  which is a light generalization of Definition \ref{def-strongly-co}. 

\begin{defn} \label{def-elem-hopf} Let $G$ be a group and let $\bar a$ be a generating tuple of $G$. We say that $G$ is \textit{elementary co-hopfian},  if there exists a formula $\varphi(\bar x)$ such that  $G \models \varphi(\bar a)$ and such that  for any endomorphism $h$ of $G$,  if $G \models \varphi(h(\bar a))$ then $h$ is an automorphism. 
\end{defn}

We emphasize  that the above definition is independent of the chosen generating tuple $\bar a$ and that a strongly co-hopfian group is elementary co-hopfian. 

\begin{lem} \label{main-lem3} Let $\G$ be a group which is either equationally noetherian or finitely presented. If $\Gamma$ is elementary  co-hopfian then $\G$ is a prime model. Moreover, if $\G$ is strongly co-hopfian then it is $\exists$-homogeneous and $\G$-determined. 
\end{lem}

\proof Let 
$$
\Gamma=\<\bar x | w_i(\bar x )=1, i \in \mathbb N\>
$$

be a presentation of $\Gamma$. Since $\Gamma$ is either equationally noetherian or finitely presented,  there exists $p \in \mathbb N $ such that 
$$
\G \models \forall \bar x ( w_1(\bar x)=1 \wedge \dots \wedge w_p(\bar x)=1 \Rightarrow w_i(\bar x)=1), \leqno (1)
$$
for any $i \in \mathbb N$.  

Let $\varphi(\bar x)$ as in Definition \ref{def-elem-hopf}. Using  $(1)$,  we conclude that for any tuple  $\bar y$ in $\G$ which satisfies 
$$
w_1(\bar y)=1 \wedge \dots \wedge w_p(\bar y)=1 \wedge \varphi(\bar y), 
$$
there is an automorphism  $f$ of $\G$ which sends $\bar x$ to $\bar y$.

Let  $\bar b$ in $\G^m$ and let us show that the orbit of $\bar b$ under the action of $Aut(\G)$ is definable. We conclude by  Proposition \ref{prop-prime-model}. 

There exists a tuples of words $\bar t(\bar x)$ such that $\bar b= \bar t(\bar x)$.  We see that the orbit of $\bar b$ is definable by the formula
$$
\psi(\bar z):= \exists \bar y (\bigwedge_{1\leq i \leq p }w_i(\bar y)=1 \wedge \varphi(\bar y) \wedge \bar z =\bar t(\bar y)).
$$

When $\G$ is strongly co-hopfian, the proof of the fact that $\G$ is $\exists$-homogenous and $\G$-determined proceed in a similar way to that of the proof of Lemma \ref{main-lem2}. \qed

\begin{rem} We note that  as a consequence,  if a finitely presented simple group is co-hopfian then it is prime.   Indeed,  if $\G$ is a finitely presented infinite simple co-hopfian group,  then by taking, in Definition \ref{def-strongly-co},  $S$ to be reduced to a one nontrivial element $g$,   then any homomorphism $\varphi : \G  \rightarrow \G$ with $\varphi(g) \neq 1$ is an automorphism and thus $\G$ is strongly co-hopfian. 
\end{rem}

We conclude this section with the following lemma of independent interest. 

\begin{lem} \label{lem-auto-mono} Let $F$ be a nonabelian free group of finite rank and let  $\bar a, \bar b$ be tuples from $F$. Then the existence of an  automorphism  sending $\bar a$ to $\bar b$  is equivalent to the existence of  a monomorphism sending $\bar a$ to $\bar b$ and   a monomorphism sending $\bar b$ to $\bar a$. 
\end{lem}

\proof  Let  $f$ and $g$ be  monomorphisms such that $f(\bar a)=\bar b$ and $g(\bar b)=\bar a$.

Then $g \circ f$ is a monomorphim which fixes $\bar a$. If $g \circ f$ is an automorphism,  then $g$ is an automorphism and we are done. If $g \circ f$ is not an automorphism then $\bar a$ is in a proper free factor of $F$ by Proposition \ref{prop1}. A similar argument can be used for $\bar b$.

So we suppose that $\bar a$ and $\bar b$ are in  proper free factors. Let $F=F_1 *A=F_2*B$ with $\bar a \in F_1$ and $\bar b \in F_2$ and such that $F_1$ (resp. $F_2$) is without proper free factor containing $\bar a$ (resp. $\bar b$).

By applying Grushko theorem to the subgroup $f(F_1)$ with respect to the decomposition $F=F_2*B$ and since $ \bar b \in f(F_1) \cap F_2$,  we get  $f(F_1)=f(F_1) \cap F_2 *K$ for some subgroup $K$ of $F$.  

We claim that $K=1$. Suppose  for  a contradiction that $K \neq 1$. Hence, by  \cite[Theorem 1.8. CH IV]{LyndonSchupp77}, $F_1$ has a decomposition $P*L$  such that  $f(P)=f(F_1) \cap F_2$ and $f(L)=K$. Since $K \neq 1$, we get  $L \neq 1$. Since $f$ is a monomorphism,  we get   $\bar a \in P$; which is clearly a contradiction with the choice of $F_1$. Thus $K=1$ and $f(F_1) \leq F_2$.

With a similar argument,  we have $g(F_2) \leq F_1$. As before,  $(g \circ f)_{|F_1}$ is a monomorphism of $F_1$ which fixes $\bar a$. If $F_1$ is cyclic then $(g \circ f)_{|F_1}$ is the identity and thus an automorphism of $F_1$. If $F_1$ is noncyclic, then, since $F_1$ is freely indecomposable relative to the subgroup generated by $\bar a$, $(g \circ f)_{|F_1}$  is an automorphism,  by Proposition \ref{prop1}. Hence $g_{|F_2}$ is surjective. In particular $F_1$ and $F_2$ have the same rank. Therefore  $f_{|F_1}$ can be extended to an automorphism of $F$. \qed

\begin{rem} Remark that Lemma \ref{lem-auto-mono} has as a consequence the $\exists$-homogeneity of the free group of rank $2$. Let $\bar u $ and $\bar v$ be tuples from $F_2$ such that $tp_{\exists}(\bar u)=tp_{\exists}(\bar v)$.  Write $F_2=\<x_1,x_2|\>$,  $\bar u=\bar w(x_1,x_2)$. Since $tp(\bar u)\subseteq tp(\bar v)$,  there exists $y_1, y_2$ such that $[y_1, y_2] \neq 1$,  $\bar v=\bar w(y_1, y_2)$. Therefore the map defined by $f(x_i)=y_i$ is a monomorphism which  sends $\bar u$ to $\bar v$. Similarly, there exists a monomorphism $g$ wich sends $\bar v$ to $\bar u$.  We conclude by Lemma \ref{lem-auto-mono}. 
\end{rem}

\section{The existential case}

We begin in this section by examining   existential types in free groups. The main purpose  is to give   the proof of Theorem \ref{thm2} with the hypothesis $P=\emptyset$.   

Let $F_1$ and $F_2$ be nonabelian free groups of finite rank and let $\bar a$ (resp. $\bar b$) be a tuple from $F_1$ (resp. $F_2$).   We denote by $Hom(F_1|\bar a,F_2|\bar b)$,  the set of homomorphisms $f : F_1 \rightarrow F_2$ such that $f(\bar a)=\bar b$.  We denote by $rk(H)$ the rank of $H$.

If $\bar a=(a_1, \dots, a_n)$ is a tuple from $F$,  we denote by $|\bar a|$ the integer 
$$
|\bar a|=\max \{|a_i | | 1 \leq i \leq n\},
$$ 
where $|a_i |$ denote the length of $a_i$ with respect to some fixed basis of the ambiant free group $F$.  In the rest of this section, we suppose that the tuples  which we use are finite and have the same length.  For a tuple $\bar a$ from $F$,  we denote by $tp_{\exists}^{F}(\bar a)$ its exiential type and by $tp_{\forall}^{F}(\bar a)$ its universal type.

\begin{defn}  Let $F_1$ and $F_2$ be nonabelian free groups of finite rank and let $\bar a$ (resp. $\bar b$) be a tuple from $F_1$ (resp. $F_2$).  We say that $(\bar a, \bar b)$ is \textit{existentially rigid},  if there is no nontrivial free decomposition $F_1=A*B$ such that $A$ contains a tuple ${\bar c}$ with $tp_{\exists}^{F_1}(\bar a) \subseteq   tp_{\exists }^A(\bar c) \subseteq tp_{\exists }^{F_2}(\bar b)$. \end{defn}

\begin{rem} $\;$

$(1)$ Since $A$ is an e.c. subgroup of $F_1$,  we have $tp_{\exists }^A(\bar c)=tp_{\exists }^{F_1}(\bar c)$.

$(2)$ We note that $(\bar a, \bar a)$ is existentially rigid if and only if $tp_{\exists}^{F_1}(\bar a)$ is not realized in any free group having a  smaller rank than the rank of $F_1$.

$(3)$ If $F$ is a free group of rank $2$, then for any tuples $\bar a, \bar b$, $(\bar a, \bar b)$ is existentially rigid. 
\end{rem}

We begin by showing the following proposition. 
\begin{prop}\label{prop-key-1} Let $F_1$ and $F_2$ be nonabelian free groups of finite rank and let $\bar a$ (resp. $\bar b$) be a tuple from $F_1$ (resp. $F_2$). Suppose that $(\bar a, \bar b)$ is existentially rigid.   Let $\bar s$ be a basis  of $F_1$. Then there exists a quantifier-free  formula $\varphi(\bar x, \bar y)$,  such that $F_1 \models \varphi(\bar a, \bar s)$ and   such that for any $f\in Hom(F_1|\bar a,F_2|\bar b)$,  if $F_2 \models \varphi(\bar b, f(\bar s))$ then $f$ is an embedding. 
\end{prop}

\proof    Let $(\psi_i(\bar x, \bar y)| i \in \mathbb N)$ be an enumeration of  the quantifier-free type of $(\bar a, \bar s)$ and set   

 $$\varphi_n(\bar x, \bar y)=\bigwedge_{0 \leq i \leq n}\psi_i(\bar x, \bar y).$$

Suppose for a  contradiction that for any $n \in \mathbb N$,  there exists a non-injective  homomorphism $f_n \in Hom(F_1|\bar a,F_2|\bar b)$ such that $F_2 \models \varphi_n(\bar b, f_n(\bar s))$. 

We emphasize the following property which will be used below implicitly. For any subsequence $(f_{n_k})_{k \in \mathbb N}$ and for any $n \in \mathbb N$, $F_2 \models \varphi_n(\bar b, f_{n_{k}}(\bar s))$ for all but finitely many $k$.

Since $f_n \in Hom(F_1|\bar a,F_2|\bar b)$,   $\bar b \in f_n(F_1)$ and since $f_n$ is not injective we get $rk(f_n(F_1))<rk(F_1)$,  for all $n$.   Using the pigeon hole principale,  we extract a subsequence, that we assume to simplify notation to be $(f_n)_{n \in \mathbb N}$ itself, such that $rk(f_n(F_1))$ is a fixed natural number $r$ for all $n$.

By Proposition \ref{cor2}, each $f_n(F_1)$ has a basis $\{d_{1n}, \dots, d_{{p_n}n}, \dots, d_{rn}\}$ such that $\bar b$ is contained in the subgroup generated by $\{d_{1n}, \dots, d_{{p_n}n}\}$ and $|d_{in}| \leq |\bar b|$ for all $1\leq i \leq p_n$ and for all $n$. 

Therefore for any $ n \in \mathbb N$, the set $\{d_{1n}, \dots, d_{p_n}\}$ is contained in the ball of radius $|\bar b|$ of $F_2$. Thus again,  using the pigeon hole principale,  we can find a subsequence, that we assume to simplify notation to be $(f_n)_{n \in \mathbb N}$ itself, such that $p_n$ is a fixed integer $p$ and $d_{in}=d_i$ for all $n \in \mathbb N$ and $1 \leq i \leq p$. 

We conclude that for all $n \in \mathbb N$,   
$$
f_n(F_1)=\<d_1, \dots, d_p, d_{(p+1)n}, \dots , p_{rn}\>
$$
and $\bar b$ is in the subgroup with basis $\{d_1, \dots, d_p\}$. 

Set $L=\<d_1, \dots, d_p, d_{(p+1)0}, \dots, d_{r0}\>=f_0(F_1)$.

\smallskip
\noindent \textit {Claim 1. $tp_{\exists }^{F_1} (\bar a) \subseteq tp_{\exists}^L(\bar b)$.}

\proof  Let $
\varphi(\bar x, \bar y)$
be a quantifier-free formula such that $F_1 \models \exists \bar y \varphi(\bar a, \bar y)$.  Then there exists a tuples of words $\bar \alpha(\bar t)$  such that $F_1 \models \varphi(\bar a, \bar \alpha(\bar s))$. By construction of the sequence $(f_n)_{n \in \mathbb N}$,   $F_2 \models \varphi(\bar b, \bar \alpha(f_n(\bar s)))$ for all but finitely many $n$. 

Since  $f_n(\bar s), \bar b $ are in $f_n(F_1)$ we get $f_n(F_1) \models \varphi(\bar b, \bar \alpha(f_n(\bar s)))$ for all but finitely many $n$.  Therefore $f_n(F_1) \models \exists \bar y \varphi(\bar b,  \bar y)$  for all but finitely many $n$. 

The map $h : f_n(F_1)  \rightarrow L$ defined by $h(d_i)=d_i$ for $1 \leq i \leq p$ and $h(d_{jn})=d_{j0}$ for $p+1 \leq j \leq r$ extend to an isomorphism which   fixes $\bar b$ and that we still denote by $h$. 

 Since $h$ is an isomorphism which fixes $\bar b$,  we conclude that  $L\models \exists \bar y  \varphi(\bar b, \bar y)$. Hence $tp_{\exists}^{F_1}(\bar a) \subseteq tp_{\exists}^{L}(\bar b)$  as claimed. \qed

\smallskip
Clearly  $tp_{\exists}^L(\bar b) \subseteq tp_{\exists}^{F_2}(\bar b)$. By Proposition \ref{prop-hom-gp-libre}, $F_1$ has a  free decomposition $F_1=A*B$ such that $f_0(A)=L$ and $f_0(B)=1$ and with $f_0$ is injective in restriction to $A$. Since $rk(L)<rk(F_1)$,  the precedent decomposition is nontrivial. 
 
Let $\bar c$ to be the unique tuple of $A$ such that $f_0(\bar c)=\bar b$.  Since $f_0$ is injective in restriction to $A$,  we get $tp_{\exists}^{A} (\bar c)= tp_{\exists}^L(\bar b)$. 
 
We conclude that $tp_{\exists}^{F_1} (\bar a) \subseteq  tp_{\exists}^A(\bar c) \subseteq tp_{\exists}^{F_2}(\bar b)$ and thus $(\bar a, \bar b)$ is not existentially rigid; a final contradiction. \qed

\begin{defn}  Let $F$ be a free group and let $\bar a=(a_1, \dots, a_m)$ be a tuple from $F$.  We say that $\bar a$ is \textit{a power of a primitive element},  if there exist integers $p_1, \dots, p_m$ and a primitive element $u$ such that $a_i=u^{p_i}$ for all $i$. 
\end{defn}

\begin{lem} \label{lem-exit-rigid}Let $F_1$ and $F_2$ be nonabelian free groups of finite rank and let $\bar a$ (resp. $\bar b$) be a tuple from $F_1$ (resp. $F_2$) such that $tp_{\exists}^{F_1}(\bar a)=tp_{\exists}^{F_2}(\bar b)$. Suppose that  $(\bar a, \bar b)$ is existentially rigid. Then either $rk(F_1)=2$  and $\bar a$  is  a power of a primitive element,  or   there exists an embedding $h :F_1 \rightarrow F_2$ such that $h(F_1)$ is an e.c. subgroup of $F_2$. 
\end{lem}

\proof  We suppose that the first case of the conclusion of the lemma is not satisfied. Let $\varphi_0 (\bar x, \bar y)$  be the quantifier-free formula given by Proposition  \ref{prop-key-1} applied to the tuple $(\bar a, \bar b)$.

Observe  that $(\bar a, \bar a)$ is also existentially rigid. Hence,  by Proposition  \ref{prop-key-1} applied to the tuple $(\bar a, \bar a)$,  we obtain   also a quantifier-free formula $\varphi_1 (\bar x, \bar y)$,   such that $F_1 \models \varphi_1 (\bar a, \bar s)$ and such that for any $f \in Hom(F_1|\bar a, F_1 |\bar a)$,  if $F_1 \models  \varphi_1(\bar a, f(\bar s))$ then $f$ is an embedding.

There exists a tuple of words $\bar w(\bar x)$ such that $\bar a= \bar w(\bar s)$.  Since $tp_{\exists}^{F_1}(\bar a) =tp_{\exists}^{F_2}(\bar b)$ we get 
$$
F_2 \models \varphi_0 (\bar b, \bar s') \wedge \varphi_1 (\bar b, \bar s') \wedge \bar b= \bar w(\bar s'),
$$
for some tuple $\bar s'$ from $F_2$. By Proposition  \ref{prop-key-1},  the map $\bar s \rightarrow \bar s'$ extend to an embedding that we denote by $h$. 

We claim that $h(F_1)$ is an e.c. subgroup of $F_2$.  Let $\psi(\bar x, \bar y)$ be an existential formula such that $F_2 \models \psi(\bar b, \bar s')$.  Then
$$
F_2 \models \exists \bar s' (\varphi_1 (\bar b, \bar s') \wedge \bar b= \bar w(\bar s') \wedge \psi(\bar b, \bar s')), 
$$
and since $tp_{\exists}^{F_1}(\bar a) =tp_{\exists}^{F_2}(\bar b)$ we get 
$$
F_1 \models \varphi_1 (\bar a, \bar s'') \wedge \bar a= \bar w(\bar s'') \wedge \psi(\bar a, \bar s''),
$$
for some tuple $\bar s''$ of $F_1$. 

Hence the map $\bar s \rightarrow \bar s''$ extend to a monomorphism of $F_1$  fixing $\bar a$ that we denote by $h'$.

By Proposition \ref{prop-hom-gp-libre}, if $h'$ is not an automorphism then $F_1$ is freely decomposable with respect to the subgroup generated by $\bar a$.  
Let $F_1=C*D$ be a nontrivial free decomposition, with $rk(C)$ of minimal rank such that $\bar a$ is in $C$. If $C$ is  nonabelian  then we get  a contradiction to the fact that $(\bar a, \bar a)$ is existentially rigid. 

Hence $C$ is abelian and in this case $C$ is cyclic.  Therefore  $\bar a$ is a power of a primitive element. We observe that if $rk(F_1)>2$,  then $F_1$ has a nonabelian free factor containing $C$ and thus $(\bar a, \bar a)$ is not existentially rigid. Therfeore   $rk(F_1)=2$; a contradiction with our assumption. 

Thus $h'$ is an automorphism of $F_1$ which fixes $\bar a$.  Therefore $F_1 \models  \psi(\bar a, \bar s)$.  Since $h$ is an embedding,  we get $h(F_1) \models \psi(\bar b, \bar s')$. Therefore $h(F_1)$ is e.c. in $F_2$ as required.  \qed

\begin{prop}  \label{prop-dicho-exit} Let $F_1$ and $F_2$ be nonabelian free groups of finite rank and let  $\bar a$ (resp. $\bar b$) be a tuple from $F_1$ (resp. $F_2$) such that $tp_{\exists}^{F_1}(\bar a)=tp_{\exists}^{F_2}(\bar b)$. Then one of the following cases holds:

\smallskip

$(1)$ There exists a tuple $\bar c$ in $F_1$ which is a power of a primitive element such that $tp_{\exists}^{F_1}(\bar a) =tp_{\exists}^{F_1}(\bar c)$; 

$(2)$ There exists an e.c. subgroup $E(\bar a)$ (resp. $E(\bar b)$) containing $\bar a$  (resp. $\bar b$) of $F_1$ (resp. $F_2$) and an isomorphism $\tau : E(\bar a) \rightarrow E(\bar b)$ sending $\bar a$ to $\bar b$. 
\end{prop}

\proof

If  $(\bar a, \bar b)$ is existentially rigid then the result follows from Lemma \ref{lem-exit-rigid}. 

Let us now treat the case when $(\bar a, \bar b)$ is not existentially rigid.  Let $F_1=C*B$ be a nontrivial free decomposition and $\bar c $ in $C$ such that $tp_{\exists}^{C}(\bar c)=tp_{\exists}^{F_1}(\bar a)$.  We may choose $C$ of minimal  rank satisfying the precedent property. 

Suppose that $C$ is freely decomposable with respect to the subgroup generated by $\bar c$. Let $C=C_1*C_2$,  with $\bar c$ is in  $C_1$.  If $C_1$ is nonabelian then $tp_{\exists}^{C_1}(\bar c)=tp_{\exists}^{C}(\bar c)$ because $C_1 \preceq_{\exists} C$.  Thus we have a contradiction with the choice of $C$ as $C_1$ has a smaller rank. 

Thus  $C_1$ is  cyclic and thus $\bar a$ has the same existential type as a power of a primitive element and we get $(1)$. 

Hence,  we assume that $C$ is freely indecomposable with respect to the subgroup generated by $\bar c$.  We see  that $(\bar c, \bar a)$ is existentially rigid in $C$ as otherwise we get a contradiction to the minimality of the rank of $C$. 

By Lemma \ref{lem-exit-rigid},  there exists an embedding $h_1 : C \rightarrow F_1$  such that $h_1(C)$ is an e.c. subgroup of $F_1$ and  $h_1(\bar c)=\bar a$. 

Similarly $(\bar c, \bar b)$ is existentially rigid and by Lemma \ref{lem-exit-rigid} there exists an embedding $h_2 : C \rightarrow F_2$  such that $h_2(C)$ is an e.c. subgroup of $F_2$ and  $h_2(\bar c)=\bar b$. 

By setting $E(\bar a)=h_1(C)$ and $E(\bar b)=h_2(C)$,   $h_2 \circ h_1^{-1} : E(\bar a) \rightarrow E(\bar b)$ is an isomorphism with $h_2 \circ h_1^{-1}(\bar a)=\bar b$ and thus we get $(2)$. \qed

\begin{prop} Let $F$ be a nonabelian free group of finite rank and let $\bar a$  and $\bar b$ be tuples from $F$  such that $tp_{\exists}^{F}(\bar a)=tp_{\exists}^{F}(\bar b)$. If $(\bar a, \bar b)$ is existentially rigid then there exists an automorphism of $F$ sending $\bar a$ to $\bar b$. 
\end{prop}

\proof By Lemma \ref{lem-exit-rigid} we treat two cases.  If  $rk(F)=2$ and $\bar a$ is a power of a primitive element,  then the result follows from the $\exists$-homogeneity of the free group of rank $2$.  The case $rk(F)=2$ and $\bar b$ is a power of a primitive element is similar. 

By Lemma \ref{lem-exit-rigid}, there exists a monomorphism sending $\bar a$ to $\bar b$ and a monomorphism sending $\bar b$ to $\bar a$. Hence we conclude by Lemma \ref{lem-auto-mono}. 
\qed

\begin{rem} We note that in the free group of rank $2$ any tuple $(\bar a, \bar b)$ is existentially rigid. Hence the above proposition can be seen as a generalisation of the $\exists$-homogeneity of the free group of rank $2$. 
\end{rem}

We need the following lemma in the proof of the next proposition. For the definition of Nielsen transformations we refer the reader to \cite{LyndonSchupp77}.  

\begin{lem}\label{e.c.subgroups} If $E$ is an e.c. subgroup of a free group of finite rank $F$ then $rk(E)\leq rk(F)$ and if $E$ is proper then $rk(E)<rk(F)$. \qed
\end{lem}

\proof  We  first claim that $E$ has a finite rank. Suppose for a contradiction that $E$ has an infinite rank and let $\{x_i | i \in \mathbb N\}$ be a basis of $E$. Let $m$ be the rank of $F$. Since $E$ is e.c. in $F$, we conclude that for every $n$ the subgroup $L_n$ generated by $\{x_1, \dots, x_n\}$ is contained in a subgroup $K_n$ of $E$ of rank at most $m$.  But each $L_n$ is also a free factor of $K_n$; which is a contradiction for big $n$. 

Hence $E$ has a finite rank $m'$. Now, as before, $E$ is contained in a subgroup of itself of rank at most $m$. Hence $m' \leq m$ as required.  

Suppose now that $E$ is proper and suppose for a contradiction that $rk(E)=rk(F)$.  Let $\{h_1, \dots, h_m\}$ be a basis of $E$ and let $\{x_1, \dots, x_m\}$ be a basis of $F$.  Then for every $i$, there exists a reduced word $w_i(\bar x)$ such that $h_i=w_i(\bar x)$.  Hence in $E$, we can find $x_1', \dots , x_m'$ such that $h_i=w_i(\bar x')$. In particular, $\{x_1', \dots, x_m'\}$ is a basis of $E$. Hence, since $\{h_1, \dots, h_m\}$ is  a basis of $E$,  there exists a sequence of Nielsen transformations sending $\{x_1', \dots, x_m'\}$ to $\{w_1(\bar x'), \dots, w_m(\bar x')\}$. The corresponding sequence of Nielsen  transformations sends $\{x_1, \dots, x_m\}$ to $\{w_1(\bar x), \dots, w_m(\bar x)\}$ in $F$,  and thus $F$ is also generated by $\{h_1, \dots, h_m\}$; a contradiction.  \qed

\begin{prop}  Let $F$ be a nonabelian free group of finite rank. Then the following properties are equivalent: 

$(1)$ $F$ is $\exists$-homogeneous; 

$(2)$ The following proerties are satisfied: 

$(i)$ If a tuple $\bar a$ is a power of a primitive element and $\bar b$ have the same existential type as  $\bar a$,  then $\bar b$ is a power of a primitive element; 

$(ii)$ Every e.c. subgroup of $F$ is a free factor. 
\end{prop}

\proof  Suppose that $(1)$ holds. We see that $(i)$ is an immediate consequence.  Let $E$ be an e.c. subgroup of $F$.  Let $\{s_1, \dots, s_p\}$ be a basis of $E$ and $\{d_1, \dots, d_q\}$ be a basis of $F$. Then by Lemma  \ref{e.c.subgroups}, $rk(E) \leq rk(F)$ and thus $p\leq q$. Let $H$ be the subgroup generated by $\{d_1, \dots, d_p\}$. Then $H$ is an e.c. subgroup of $F$ and thus $tp_{\exists}^{F}(d_1, \dots, d_q)=tp_{\exists}^{F}(s_1, \dots, s_q)$.  Hence by $(1)$,  there is an automorphism sending $E$ to $H$ and thus $E$ is a free factor.

Suppose that $(2)$-$(i)$-$(ii)$ hold.  The case of powers of primitive elements is resolved by $(i)$ and the other case is resolved by $(ii)$ using Proposition \ref{prop-dicho-exit}. 
\qed

\section{Homogeneity in free groups}

We are concerned in this section with homogeneity  in free groups and  the main purpose is to give the proof of Theorem \ref{thm1} with the  hypothesis $P=\emptyset$. The general case will be treated in the next section. We use  notation of the precedent section.  For a tuple $\bar a$ from $F$,  we denote by $tp_{\exists \forall}^{F}(\bar a)$ its $\exists \forall$-type

\begin{defn}  Let $F$   be nonabelian free group of finite rank and let $\bar a$  be a tuple of $F$.  We say that $\bar a$ is \textit{rigid} if there is no nontrivial free decomposition $F=A*B$ such that $A$ contains a tuple ${\bar c}$ with $tp_{\exists \forall}^{F_1}(\bar a) =   tp_{\exists \forall }^A(\bar c)$. \end{defn}

The first purpose is to show  the following proposition,  which is the analogue of  Proposition \ref{prop-key-1}. But before it,  we shall need a preliminary study of certain sequences of subgroups similar to those who appear in the proof of  Proposition \ref{prop-key-1}. 

\begin{prop}\label{prop-key2} Let $F_1$ and $F_2$ be nonabelian free groups of finite rank and let $\bar a$ (resp. $\bar b$) be  a tuple from $F_1$ (resp. $F_2$) such that $tp_{\exists \forall}^{F_1}(\bar a)=tp_{\exists \forall}^{F_2}(\bar b)$. Suppose that  $\bar a$ is rigid in $F_1$ and let $\bar s$ be a basis of $F_1$. Then there exists an universal    formula $\varphi(\bar x, \bar y)$ such that $F_1 \models \varphi(\bar a, \bar s)$ and   such that for any $f\in Hom(F_1|\bar a,F_2|\bar b)$,  if $F_2 \models \varphi(\bar b, f(\bar s))$ then $f$ is an embedding. 
\end{prop}

\begin{defn}  \label{def-bad-sequence} Let $F$ be a free group and let  $\bar b$ be a tuple from $F$. A sequence    $(L_n| n \in \mathbb N)$ of subgroups of $F$  is called \textit{good} if it satisfies the following properties:  

\smallskip
$(1)$ There exists a fixed group $D$ such that:

$\;$$\;$$\;$$\;$$(i)$ $D$ contains  $\bar b$;

$\;$$\;$$\;$$\;$$(ii)$ $D$ is freely indecomposable relative to the subgroup generated by $\bar b$; 

$\;$$\;$$\;$$\;$$(iii)$ $D$ is a free factor of $L_n$ for all $n$;

$(2)$ There exists a fixed integer $r$ such $rk(L_n) = r$ for all $n$;

$(3)$ For any universal formula $\varphi(\bar x, \bar y) $ such that  $\exists \bar y \varphi(\bar x, \bar y)\in tp_{\exists \forall}^{F}(\bar b)$,  there exists  $n \in \mathbb N$ and $\bar \alpha_n  \in L_n$ such that   $F\models  \varphi(\bar b, \bar \alpha_n)$. 

\smallskip

For such a sequence, $r$ is called the \textit{rank} and  $D$ is called the \textit{free factor}. 
\end{defn}

Our aim now is to show the following proposition. 

\begin{prop} \label{prop-good-sequences}Let $F$ be a nonabelian free group of finite rank and let $\bar b$ be a tuple from $F$. If $(L_n| n \in \mathbb N)$ is a good sequence then there exists $p$ and a tuple $\bar c$ from $L_p$ such that $tp_{\exists \forall}^F(\bar b)=tp_{\exists \forall}^{L_p}(\bar c)$. 

\end{prop}

Before proving the previous proposition, we shall need a preliminary work on  properties of  good  sequences and powers of primitive elements.

\begin{lem} \label{lem-forall-exists}Let $F$ be a nonabelian free group of finite rank and let $\bar b$ be a tuple from $F$.   If $(L_n| n \in \mathbb N)$ is a good sequence then $tp_{\exists \forall}^F(\bar b) \subseteq tp_{\exists \forall}^{L_n}(\bar b)$ for all $n$. 
\end{lem}

\proof Let $
\varphi(\bar x, \bar y, \bar z)$
be a quantifier-free formula such that $F \models \exists \bar y \forall \bar z\varphi(\bar b, \bar y, \bar z)$. 

By Definition \ref{def-bad-sequence}(3),   we have $F \models \forall \bar z\varphi(\bar b, \bar \beta, \bar z)$  for some  $p$ and a tuple $\bar \beta$ in $L_p$. 

Since the precedent formula is universal and $\bar \beta, \bar b $ are in $L_p$,  we obtain $L_p \models \forall \bar z\varphi(\bar b, \bar \alpha, \bar z)$.  Therefore $L_p \models \exists \bar y \forall \bar z\varphi(\bar b,  \bar y, \bar z)$. 

By Definition \ref{def-bad-sequence}(1), $L_n=D*C_n$ for all $n$ and by Definition \ref{def-bad-sequence}(2) we have $rk(C_n)=rk(C_m)$ for all $n,m$. 

Therefore for any $n$,  there exists an isomorphism $h_n : L_n \rightarrow L_p$ fixing $D$ pointwise. Since $h_n$ is an isomorphism fixing $\bar b$,  we get for all $n$, $L_n \models  \exists \bar y \forall \bar z\varphi(\bar b,  \bar y, \bar z)$  as required. \qed

\begin{lem} \label{lem-power-primitive}Let $F$ be a nonabelian free group of finite rank. If $\bar a$ is a power of a primitive element and $\bar b$ is such that $tp_{\exists}^F(\bar a)=tp_{\exists}^F(\bar b)$ then $tp_{\exists \forall}^F(\bar a)\subseteq tp_{\exists \forall }^F(\bar b)$. 
\end{lem}

\proof 

Write $\bar a=(a_1, \dots, a_q)$ and $\bar b=(b_1, \dots, b_q)$. First we prove

\bigskip
\noindent
\textit{Claim 1. There exists a primitive element $u$ and an element $v$ such that: }

\smallskip
\textit{$(i)$ $tp_{\exists}^F(v)=tp_{\exists}^F(u)$; }

\smallskip
\textit{$(ii)$ There are integers $p_1, \dots, p_q$ such that $a_i=u^{p_i}$ and $b_i=v^{p_i}$ for all $i$. }

\proof

Let $u$ be a primitive element and let $p_1, \dots, p_q$ such that $a_i=u^{p_i}$ for all $i$.  Since $tp_{\exists}^F(\bar a)=tp_{\exists}^F(\bar b)$,  we find $v \in F$ such that $b_i=v^{p_i}$ for all $i$. Let $\varphi(x) \in tp_{\exists}^F(u)$.  Then $$F \models \exists v' (\varphi (v') \wedge_{1\leq i \leq q} b_i=v'^{p_i}),$$
and since $F$ is torsion-free and commutative transitive,  we conclude that $v=v'$ and thus $\varphi(x) \in tp_{\exists}^F(v)$. The inclusion $tp_{\exists}^F(v) \subseteq tp_{\exists}^F(u)$ can be proved using a similar argument. 
\qed

Let $\{x_1, \dots, x_n\}$ be a basis of $F$ and let $L$ to be the free group with basis $\{x_1, \dots, x_n, d\}$. Now we show the following claim.

\bigskip
\noindent
\textit{Claim 2.  Let $u$ and $v$ as in Claim 1.  Then $ tp_{\exists}^L(u,d)= tp_{\exists}^L(v,d)$.  }

\proof

First we show that $ tp_{\exists}^L(u,d)\subseteq tp_{\exists}^L(v,d)$.  Let us denote by $E(u,d)$ (resp. $E(v,d)$) the subgroup generated by $\{u,d\}$ (resp. $\{v,d\}$). 

Let $\varphi(x,y) \in tp_{\exists}^L(u,d)$. Since $E(u,d)$ is a free factor of $L$ it is an e.c. subgroup of $L$ and thus $E(u,d) \models \varphi(u,d)$.  Since $E(u,d)$ and $E(v,d)$ are isomorphic by the map sending $u$ to $v$ and $d$ to itself,  we conclude that $E(v,d) \models \varphi(v,d)$ and therefore $L \models \varphi(v,d)$ as required.

Now we show that $ tp_{\exists}^L(v,d)\subseteq tp_{\exists}^L(u,d)$. Let $\varphi(x,y) \in tp_{\exists}^L(v,d)$.  Then $\varphi(x,y)$ can be written as 
$$
\exists \bar z \bigvee_{1 \leq i \leq p}(\bigwedge_{W \in P_i}W(x,y, \bar z)=1\wedge \bigwedge_{V \in N_i}V(x,y,\bar z) \neq 1),
$$
where $P_i$ and $N_i$ are finite  for all $i$.  Hence there is a tuple of words $\bar \alpha(\bar x, t)$ and $p$ such that
$$
L \models \bigwedge_{W \in P_p}W(v,d, \bar \alpha(\bar x, d))=1\wedge \bigwedge_{V \in N_p}V(v,d,\bar \alpha (\bar x, d)) \neq 1. 
$$

Now we have the following observation. Let $v(\bar x)$ be a reduced word such that $v=v(\bar x)$ in $F$. Then $L$ can be viewed as the group with the generating set $\{x_1, \dots, x_n, d,v\}$ and with the presentation $v=v(\bar x)$.  Hence in any group $G$ with  a generating set $\{x_1', \dots, x_n', d_0, v_0\}$ satisfying $v_0=v(\bar x')$ we get
$$
\bigwedge_{W \in P_p}W(v_0,d_0, \bar \alpha(\bar x', d_0))=1. 
$$

Since $F$ is an e.c. subgroup of $L$  and since $tp_{\exists}^F(v)=tp_{\exists}^F(u)$, we find $\bar x', d' \in F$ such that 
$$
F \models \bigwedge_{W \in P_p}W(u,d', \bar \alpha(\bar x', d'))=1\wedge \bigwedge_{V \in N_p}V(u,d',\bar \alpha (\bar x', d')) \neq 1 \wedge u=v(\bar x').
$$
Let $G$ be the subgroup of $L$ generated by $\{x_1', \dots, x_n',   d, u\}$. Since $u=v(\bar x')$ we get by the above observation and by replacing $v_0$ by $u$ and $d_0$ by $d$ that  
$$
L \models \bigwedge_{W \in P_p} W(u,d, \bar \alpha(\bar x', d))=1. \leqno (1)
$$

Let $f :L \rightarrow F$ be the homomorphism fixing pointwise $F$ and sending $d$ to $d'$. Since 
$$
F \models  \bigwedge_{V \in N_p}V(u,f(d),\bar \alpha (\bar x', f(d))) \neq 1,
$$
we conclude that 
$$
L \models  \bigwedge_{V \in N_p}V(u,d,\bar \alpha (\bar x', d)) \neq 1. \leqno (2)
$$

By $(1)$ and $(2)$ we conclude that 
$$
L \models \bigwedge_{W \in P_p} W(u,d, \bar \alpha(\bar x', d))=1 \wedge \bigwedge_{V \in N_p}V(u,d,\bar \alpha (\bar x', d)) \neq 1,
$$
and finally $L \models \exists \bar z (\bigwedge_{W \in P_p} W(u,d, \bar z)=1 \wedge \bigwedge_{V \in N_p}V(u,d,\bar z) \neq 1).$ 

Thus $ tp_{\exists}^L(v,d)\subseteq tp_{\exists}^L(u,d)$ as required and this ends the proof of the claim. 
\qed

\bigskip
\noindent
\textit{Claim 3.  Let $u$ and $v$ as in Claim 1.  Then $ tp_{\exists \forall}^L(u,d) \subseteq  tp_{\exists \forall}^L(v,d)$.  }

\proof Let us denote by $E(u,d)$ (resp. $E(v,d)$) the subgroup generated by $\{u,d\}$ (resp. $\{v,d\}$).  By Claim 1, $ tp_{\exists}^L(u,d)= tp_{\exists}^L(v,d)$ and since $E(u,d)$ is en e.c. subgroup of $L$ we conclude that $E(v,d)$ is e.c. in $L$. 

Since $E(u,c)$ is an elemenatry subgroup of $L$ we have $ tp_{\exists \forall}^L(u,d)=  tp_{\exists \forall}^{E(u,d)}(u,d)$.  Since $E(u,d)$ and $E(v,d)$ are isomorphic by the map sending $u$ to $v$ and fixing $d$ we conclude that $ tp_{\exists \forall}^{E(u,d)}(u,d)=  tp_{\exists \forall}^{E(v,d)}(v,d)$.  Therefore $ tp_{\exists \forall}^{L}(u,d)=  tp_{\exists \forall}^{E(v,d)}(v,d)$.

Now since $E(v,d)$ is e.c. in $L$ we get  $tp_{\exists \forall}^{E(v,d)}(v,d) \subseteq  tp_{\exists \forall}^{L}(v,d)$ and finally $tp_{\exists \forall}^{L}(u,d) \subseteq  tp_{\exists \forall}^{L}(v,d)$ as required. \qed

\smallskip
Let $u$ and $v$ as in Claim 1. It follows by Claim 3 that $tp_{\exists \forall}^{L}(u) \subseteq  tp_{\exists \forall}^{L}(v)$.  Since $F$ is an elementary subgroup of $L$ we conclude that $tp_{\exists \forall}^{F}(u) \subseteq  tp_{\exists \forall}^{F}(v)$. 

Now let us show that $tp_{\exists \forall}^{F}(\bar a) \subseteq  tp_{\exists \forall}^{F}(\bar b)$. Let $p_1, \dots, p_q$ given by Claim 1. Let $\varphi(x_1, \dots, x_q) \in tp_{\exists \forall}^{F}(\bar a) $.  Then $\varphi(x^{p_1}, \dots, x^{p_q}) \in tp_{\exists \forall }^F(u)$. Since $tp_{\exists \forall}^{F}(u) \subseteq  tp_{\exists \forall}^{F}(v)$ we conclude that $\varphi(x^{p_1}, \dots, x^{p_q}) \in tp_{\exists \forall }^F(v)$ and hence  $\varphi(x_1, \dots, x_q) \in tp_{\exists \forall}^{F}(\bar b)$. This ends the proof of the lemma. \qed

\smallskip
Having disposed of this preliminary step, we are now in a position  to prove Proposition \ref{prop-good-sequences}. 

\bigskip
\noindent
\textbf{Proof of  Proposition \ref{prop-good-sequences}}. 

The proof is by induction on the rank of good sequences. Let $(L_n| n \in \mathbb N)$ be a good  sequence of $F$ and let $r$ be its rank and let $D$ be its free factor. Let $\bar s$ be a basis of $F$.     Let $(*)$  be the following property:

\smallskip
$(*)$ \textit{for any universal formula $\varphi (\bar x)$ such that $F \models \varphi(\bar s)$,  there exists $f \in Hom(F|\bar b, F|\bar b)$ such that 
$F \models \varphi(f(\bar s))$ with $f$ is non-injective in restriction to $D$. }

\smallskip
We are going to handle two cases according to $(*)$ is or not satisfied.  Let us first treat the case when $(*)$ holds. 

\bigskip
\noindent
\textit{Claim 1. There exists a sequence $(H_p|p \in \mathbb N)$ satisfying the following properties: }

\smallskip
\textit{$(i)$ For any $p \in \mathbb N$,  there exist $n \in \mathbb N$ and   $f \in Hom(F|\bar b, F|\bar b)$ such that $H_p=f(L_n)$ and such that $f$ is non-injective in restriction to $L_n$; }

\smallskip
\textit{$(ii)$ For any universal formula $\varphi(\bar x, \bar y)$ such that $\exists \bar y  \varphi(\bar x, \bar y)\in tp_{\forall \exists}^F(\bar b)$,   there exists $p_0$ such that for any $p \geq p_0$ there exists $\bar \beta_p$ in $H_p$ such that $F \models \varphi(\bar b, \bar \beta_p)$. }

\proof

Let $(\psi_i(\bar x, \bar y_i)| i \in \mathbb N)$ be an enumeration of universal formula such that $\exists \bar y_i \varphi(\bar x, \bar y_i) \in tp_{\forall}^{F}(\bar s)$ and let   for every $n \in \mathbb N$, 

 $$\varphi_n(\bar x, \bar y_0, \dots, \bar y_n)=\bigwedge_{0 \leq i \leq n}\psi_i(\bar x, \bar y_i).$$
 
 We define  $(H_p|p \in \mathbb N)$ as follows. Let $p \in \mathbb N$. Since $(L_n | n \in \mathbb N)$ is good, by $(3)$ of Definition \ref{def-bad-sequence},  there exists $n_p \in \mathbb N$ such that for some sequence $(\bar \alpha_0, \dots, \bar \alpha_p)$ in $L_{n_p}$,  
 $$
 F \models \varphi_p(\bar b, \bar \alpha_0, \dots, \bar \alpha_p). 
 $$

By $(*)$,   there exists a   homomorphism $f \in Hom(F|\bar b,F|\bar b)$ such that $$F \models \varphi_p(\bar b, f(\bar \alpha_0), \dots, f(\bar \alpha_p))$$ which is not injective in restriction to $D$. In particular $f$ is non-injective in restriction to $L_{n_p}$. 

Put $H_p=f(L_{n_p})$. Thus we get $(i)$. We note that $(f(\bar \alpha_0), \dots, f(\bar \alpha_p))$ is a sequence of $H_p$.  By construction we have $(ii)$. \qed

\smallskip
We notice that by construction  any subsequence $(H_{p_k}| k \in \mathbb N)$ satisfies also $(i)$ and $(ii)$ of Claim 1.

By $(i)$,  $\bar b \in H_p=f(L_q)$ for some $q$, and since $f$ is not injective in restriction to $D$ we have $rk(H_p)<rk(L_q)=r$ for all $p$.   Proceeding as in the proof of Proposition \ref{prop-key-1}, by using the pigeon hole principale we extract a subsequence, that we assume to simplify notation to be $(H_p|p \in \mathbb N)$ itself, such that $rk(H_p)$ is a fixed natural number $r'<r$ for all $p$.

Again, proceeding as in the proof of Proposition \ref{prop-key-1}, and up to exctracting a subsequence, we may assume that  for all $p \in \mathbb N$,  
$$
H_p=\<h_1, \dots, h_q, h_{(q+1)p}, \dots , h_{rp}\>
$$
and $\bar b$ is in the subgroup with basis $\{h_1, \dots, h_q\}$.

Let $H$ to be the subgroup with basis $\{h_1, \dots, h_q\}$. By the Grushko decomposition,  we have $H=M*N$ where $\bar b$ is in $M$ and $M$ is freely indecomposable with respect to the subgroup generated by $\bar b$.  We define $M$ to be the free factor of the sequence and thus we get $(1)$ of Definition \ref{def-bad-sequence}.   By construction,  the sequence $(H_p | p \in \mathbb N)$ satisfies $(2)$ and $(3)$ of    Definition \ref{def-bad-sequence} and hence it is a good sequence. 

By induction,  there exists $p$ such that $H_p$ has a tuple $\bar c $ with $tp_{\exists \forall}^F(\bar b)=tp_{\exists \forall}^{H_p}(\bar c)$. Now by construction,  there exists $q$ and $f \in Hom(F|\bar b,F|\bar b)$ such that $H_p=f(L_q)$ and $f$ is non-injective in restriction to $L_q$.  By Proposition \ref{prop-hom-gp-libre}, $L_q$  has a  free decomposition $L_q=A*B$ such that $f(A)=H_p$ and $f(B)=1$ and with $f$ is injective in restriction to $A$.  Since $f$ is injective in restriction to $A$,  $A$ contains  a tuple $\bar c' $ such that $tp_{\exists \forall}^{H_p}(\bar c)=   tp_{\exists \forall }^{A}(\bar c')$. Since $A$ is an elementary subgroup of $L_p$ and  $tp_{\exists \forall}^F(\bar b)=tp_{\exists \forall}^{H_p}(\bar c)$,  we conclude that $tp_{\exists \forall}^{L_q}(\bar c')=   tp_{\exists \forall }^{F}(\bar b)$.  This ends the proof when $(*)$ is satisfied.

 We treat now  the case when $(*)$ is not true.  We treat the two cases depending on the fact that $D$ is abelian or not. 

Suppose that $D$ is abelian. Hence $D$ is cyclic and we assume that it is generated by $u$. Write $\bar b=(b_1, \dots, b_q)$ and let $p_1, \dots, p_q$ integers such that $b_i=u^{p_i}$ for all $i$. 

Let $u'$ be a primitive element in $F$ and let $\bar b'=(u'^{p_1}, \dots, u'^{p_q})$. By Theorem \ref{thm-elemantary-free-factor}, we conclude that $tp_{\exists \forall}^{L_n}(\bar b)=tp_{\exists \forall}^{F}(\bar b')$. 

 By Lemma \ref{lem-forall-exists}, $tp_{\exists \forall}^F(\bar b) \subseteq tp_{\exists \forall}^{L_n}(\bar b)$ for all $n$.  Therefore $tp_{\exists \forall}^F(\bar b) \subseteq tp_{\exists \forall}^{F}(\bar b')$. In particular $tp_{\exists }^F(\bar b) = tp_{\exists }^{F}(\bar b')$.  By Lemma \ref{lem-power-primitive},  we get $tp_{\exists \forall}^F(\bar b) = tp_{\exists \forall}^{F}(\bar b')$. 
 
 By $tp_{\exists \forall}^{L_n}(\bar b)=tp_{\exists \forall}^{F}(\bar b')$,  we conclude that $tp_{\exists \forall}^{L_n}(\bar b)=tp_{\exists \forall}^{F}(\bar b)$. Hence in this case we get the required result.

Suppose now  that $D$ is nonabelian. Since $(*)$ is not true,  there exists a universal formula $\varphi_0(\bar x)$ such that $F \models \varphi(\bar s)$ and such that for any $f \in Hom(F|\bar b, F|\bar b)$ if  
$F \models \varphi_0(f(\bar s))$  then $f$ is injective  in restriction to $D$.

We claim that $D$ is e.c. in $F$. Let $\bar d$ be a basis of  $D$. 
Then there exists a tuple of words $\bar w(\bar s)$ such that $\bar d= \bar w(\bar s)$ and a tuple of words $\bar v (\bar y)$ such that $\bar b= \bar v(\bar d)$.

 
 Let $\varphi(\bar x, \bar y , \bar z)$ be a quantifier-free formula such that $F\models \exists \bar z \varphi(\bar b, \bar d , \bar z)$. 

Thus 
$$
F \models   \exists \bar z \exists \bar d \exists \bar s (\varphi(\bar b, \bar d , \bar z)\wedge \varphi_0(\bar s) \wedge \bar d= \bar w(\bar s) \wedge \bar b= \bar v(\bar d)). 
$$

Since $(L_n| n \in \mathbb N)$ is a good sequence,  there exist $p$ and tuples of elements of $L_p$, $\bar \alpha$, $\bar d'$, $\bar s'$ such that
$$
F \models \varphi(\bar b, \bar d' , \bar \alpha)\wedge \varphi_0(\bar s') \wedge \bar d'= \bar w(\bar s') \wedge \bar b= \bar v(\bar d')). 
$$

Hence the homomorphism $f$ which sends $\bar s$ to $\bar s'$ is injective on $D$ and fixes $\bar b$.  

Let $D'$ to be the subgroup of $L_p$ generated by $\bar d'$.  Using the Grushko decomposition and since $\bar b$ is in $D \cap D'$ and since $D$ is freely indecomposbale relative to the subgroup generated by $\bar b$,  we conclude that $D' \leq D$. Therefore the  map $\bar d \rightarrow \bar b'$ extend to a monomorphism $h$ of $D$ fixing $\bar b$.  Since $D$ is freely indecomposable relative to the subgroup generated by $\bar b$, by Proposition \ref{prop1} $h$ is an automorphism of $D$. Since $D$ is a free factor of $L_n$,  $h$ can be extended to an automorphism of $L_n$ that we still denote by $h$.   

Since 
$$
L_n \models  \exists \bar z \varphi(\bar b, h(\bar d) , \bar z),
$$
we conclude that 
$$
L_n \models  \exists \bar z \varphi(\bar b, \bar d, \bar z),
$$
and thus $D \models \exists \bar z \varphi(\bar b, \bar d, \bar z)$ as $D$ is e.c. in $L_n$.

Hence $D$ is an e.c. subgroup of $F$ as claimed. Thus $tp_{\exists \forall }^D(\bar b) \subseteq tp_{\exists \forall }^F(\bar b)$ and since $D$ is e.c. subgroup of $L_n$ we get  $tp_{\exists \forall }^{L_n}(\bar b) \subseteq tp_{\exists \forall }^F(\bar b)$. Therefore  $tp_{\exists \forall }^{L_n}(\bar b) = tp_{\exists \forall }^F(\bar b)$ by Lemma \ref{lem-forall-exists}. This ends the proof of the proposition. \qed 

\smallskip
\noindent 
\textbf{Proof of Proposition \ref{prop-key2}}.

The proof proceeds in a similar way to that of Proposition \ref{prop-key-1}. Let $(\psi_i(\bar x, \bar y)| i \in \mathbb N)$ be an enumeration of $tp_{\forall}^{F_1}(\bar a, \bar s)$ and set  $$\varphi_n(\bar x, \bar y)=\wedge_{1 \leq i \leq n}\psi_i(\bar x, \bar y).$$

Suppose for a  contradiction that for any $n \in \mathbb N$,  there exists a non-injective  homomorphism $f_n \in Hom(F_1|\bar a,F_2|\bar b)$ such that $F_2 \models \varphi_n(\bar b, f_n(\bar s))$. 

Observe  that for any subsequence $(f_{n_k})_{k \in \mathbb N}$ and for any $n \in \mathbb N$,  there exists $n_k$ such that for any $k' \geq k$ we have $F_2 \models \varphi_n(\bar b, f_{n_{k'}}(\bar s))$.

We have $\bar b \in f_n(F_1)$ and since $f_n$ is not injective we have $rk(f_n(F_1))<rk(F_1)$ for all $n$.   Using the pigeon hole principale,  we extract a subsequence, that we assume to simplify notation to be $(f_n)_{n \in \mathbb N}$ itself, such that $rk(f_n(F_1))$ is a fixed natural number $r$ for all $n$.

Proceeding as in the proof of Proposition \ref{prop-key-1}, and up to exctracting a subsequence, we may assume that   for all $n \in \mathbb N$,   
$$
f_n(F_1)=\<d_1, \dots, d_p, d_{(p+1)n}, \dots , p_{rn}\>
$$
and $\bar b$ is in the subgroup with basis $\{d_1, \dots, d_p\}$. 

Set $L_n=\<d_1, \dots, d_p, d_{(p+1)n}, \dots, d_{rn}\>=f_n(F_1)$. Let $H$ to be the subgroup with basis $\{d_1, \dots, d_q\}$. By the Grishko decomposition,   $H=D*N$ where $\bar b$ is in $D$ and $D$ is freely indecomposable with respect to the subgroup generated by $\bar b$. 

We claim now that the sequence $(L_n | n \in \mathbb N)$ is a good sequence. By construction, $(L_n | n \in \mathbb N)$ satisfes $(1)$ and $(2)$ of  Definition \ref{def-bad-sequence} and it remains to show $(3)$ of the same definition. 

Let $
\varphi(\bar x, \bar y)$
be an universal  formula such that $F_2 \models \exists \bar y \varphi(\bar b, \bar y)$.  Since $tp_{\exists \forall}^{F_1}(\bar a)=tp_{\exists \forall}^{F_2}(\bar b)$ we get 
$$
F_1 \models \exists \bar y \varphi(\bar a, \bar y).
$$

Therefore there exists a tuples of words $\bar \alpha(\bar t)$  such that $F_1 \models \varphi(\bar a, \bar \alpha(\bar s))$. By construction of the sequence $(f_n)_{n \in \mathbb N}$ we have $F_2 \models \varphi(\bar b, \bar \alpha(f_n(\bar s)))$ for all but finitely many $n$. 

Therefore for a large $n$ we have a tuple  $\bar \alpha_n=\alpha(f_n(\bar s))$ in $L_n$ such that $F_2 \models \varphi(\bar b, \bar \alpha_n)$ and thus we get $(3)$ of Definition \ref{def-bad-sequence}. 

We conclude that $(L_n | n \in \mathbb N)$ is a good sequence as claimed. By Proposition \ref{prop-good-sequences}, there exists $p$ and a tuple $\bar c$ from $L_p$ such that $tp_{\exists \forall}^{F_2}(\bar b)=tp_{\exists \forall}^{L_p}(\bar c)$.

 By Proposition \ref{prop-hom-gp-libre}, $F_1$  has a  free decomposition $F_1=A*B$ such that $f_p(A)=L_p$ and $f(B)=1$ and with $f_p$ is injective in restriction to $A$.  Since $f_p$ is not injective the above decomposition is nontrivial. 

Thus $A$ has a tuple  $\bar c'$ with   $tp_{\exists \forall}^{F_2}(\bar b)=   tp_{\exists \forall }^{A}(\bar c')$. Since $A$ is an elementary subgroup of $F_1$ and since $tp_{\exists \forall}^{F_2}(\bar b)=tp_{\exists \forall}^{F_1}(\bar a)$,  we conclude that $\bar a$ is not rigid,  which is our final  contradiction. This ends the proof of the proposition.  \qed

\smallskip
The following proposition  is the analogue of Lemma \ref{lem-exit-rigid}. 

\begin{prop} \label{prop-cas-rigid} Let $F_1$ and $F_2$ be nonabelian free groups of finite rank and let $\bar a$ (resp. $\bar b$) be a tuple from $F_1$ (resp. $F_2$) such that $tp_{\exists \forall}^{F_1}(\bar a)=tp_{\exists \forall}^{F_2}(\bar b)$. Suppose that  $\bar a$ is  rigid. Then either $rk(F_1)=2$  and $\bar a$  is  a power of a primitive element,  or  there exists an embedding $h :F_1 \rightarrow F_2$ such that $h(F_1) \preceq_{\exists \forall} F_2$.
\end{prop}

\proof  We suppose that the first case of the conclusion of the proposition is not satisfied. Let $\varphi_0 (\bar x, \bar y)$  be the universal  formula given by Proposition  \ref{prop-key2} applied to the tuple $(\bar a, \bar b)$.

By Proposition  \ref{prop-key2} applied to the tuple $(\bar a, \bar a)$,  we get also a universal formula $\varphi_1 (\bar x, \bar y)$  such that $F_1 \models \varphi_1 (\bar a, \bar s)$ and such that for any $f \in Hom(F_1|\bar a, F_1 |\bar a)$ if $F_1 \models  \varphi_1(\bar a, f(\bar s))$ then $f$ is an embedding.

There exists a tuple of words $\bar w(\bar x)$ such that $\bar a= \bar w(\bar s)$.  Since $tp_{\exists \forall}^{F_1}(\bar a) =tp_{\exists \forall}^{F_2}(\bar b)$ we get 
$$
F_2 \models \varphi_0 (\bar b, \bar s') \wedge \varphi_1 (\bar b, \bar s') \wedge \bar b= \bar w(\bar s'),
$$
for some tuple $\bar s'$ in $F_2$. By Proposition  \ref{prop-key2},  the map $\bar s \rightarrow \bar s'$ extend to an embedding that we denote by $h$. 

We claim that $h(F_1)\preceq_{\exists \forall}F_2$.  Let $\psi(\bar x, \bar y)$ be $\exists \forall$-formula such that $F_2 \models \psi(\bar b, \bar s')$.  Then
$$
F_2 \models \exists \bar s' (\varphi_1 (\bar b, \bar s') \wedge \bar b= \bar w(\bar s') \wedge \psi(\bar b, \bar s')), 
$$
and since $tp_{\exists \forall}^{F_1}(\bar a) =tp_{\exists \forall}^{F_2}(\bar b)$,  we get 
$$
F_1 \models \varphi_1 (\bar a, \bar s'') \wedge \bar a= \bar w(\bar s'') \wedge \psi(\bar a, \bar s''),
$$
for some tuple $\bar s''$ of $F_1$. 

Hence the map $\bar s \rightarrow \bar s''$ extend to a monomorphism of $F_1$  fixing $\bar a$ that we denote by $h'$.

By Proposition \ref{prop-hom-gp-libre}, if $h'$ is not an automorphism then $F_1$ is freely decomposable with respect to the subgroup generated by $\bar a$.  
Let $F_1=C*D$ be a nontrivial free decomposition with $rk(C)$ of minimal rank such that $\bar a$ is in $C$. If $C$ is  nonabelian  then we get to a contradiction to the fact that $\bar a$ is  rigid. 

Hence $C$ is abelian and in this case $D$ is cyclic and thus $\bar a$ is a power of a primitive element and $rk(F_1)=2$; a contradiction with our assumption. 

Thus $h'$ is an automorphism of $F_1$ which fixes $\bar a$.  Therefore $F_1 \models  \psi(\bar a, \bar s)$. 

Since $h$ is an embedding we get $h(F_1) \models \psi(\bar b, \bar s')$. Therefore $h(F_1) \preceq_{\exists \forall}F_2$ as required.   \qed

We give now  the proof of Theorem \ref{thm1},  with the hypothesis $P=\emptyset$. 

\begin{prop} \label{prop-princip}Let $F$ be a nonabelian free group of finite rank  and let $\bar a$  and $\bar b$ be tuples of $F$ such that $tp^F(\bar a)=tp^F(\bar b)$. Then there exists an automorphism $\sigma$ of $F$ such that $\sigma(\bar a)=\bar b$. 
\end{prop}

\proof 

We may assume that $rk(F)>2$. Suppose that $\bar a$ is  rigid. It follows in particular that $\bar b$ is also rigid.  By Proposition \ref{prop-cas-rigid}, there exists a monomorphism sending $\bar a$ to $\bar b$ and a monomorphism sending $\bar b$ to $\bar a$. Hence,  we conclude by Lemma \ref{lem-auto-mono}. 

We treat now the case $\bar a$ is not   rigid.  Let $F_1=C*B$ be a nontrivial free decomposition and $\bar c $ in $C$ such that $tp_{\exists \forall}^{C}(\bar c)=tp_{\exists \forall}^{F_1}(\bar a)$.  We may choose $C$ of  minimal rank satisfying the precedent property. 

Suppose that $C$ is freely decomposable with respect to the subgroup generated by $\bar c$. Let $C=C_1*C_2$ with $\bar c$ is in  $C_1$.  If $C_1$ is nonabelian then $tp_{\exists \forall}^{C_1}(\bar c)=tp_{\exists \forall}^{C}(\bar c)$,  because $C_1 \preceq_{\exists \forall} C$ by Theorem \ref{thm-elemantary-free-factor}.  Thus we have a contradiction with the choice of $C$ as $C_1$ has a smaller rank. 

Thus  $C_1$ is  cyclic and thus $\bar a$ has the same $\exists \forall$-type as a power of a primitive element.  By Theorem \ref{thm-elemantary-free-factor}, we conclude that $\bar a$ has the same type as a power of a primitive element and by Theorem \ref{primitive},  we get the required conclusion. 

Hence,  we assume that $C$ is freely indecomposable with respect to the subgroup generated by $\bar c$.  We see  that $\bar c$ is  rigid in $C$ as otherwise we get a contradiction to the minimality of the rank of $C$. 

By Proposition \ref{prop-cas-rigid},  there exists an embedding $h_1 : C \rightarrow F_1$  such that $h_1(C) \preceq _{\exists \forall}F$ and  $h_1(\bar c)=\bar a$. 

Similarly,  by Proposition \ref{prop-cas-rigid},  there exists an embedding $h_2 : C \rightarrow F$  such that $h_2(C) \preceq_{\exists \forall}F$ and  $h_2(\bar c)=\bar b$. 

 We have $h_2 \circ h_1^{-1} : h_1(C) \rightarrow h_2(C)$ is an isomorphism with $h_2 \circ h_1^{-1}(\bar a)=\bar b$. 
 Since $h_2(C) \preceq_{\exists \forall}F$ and $h_1(C) \preceq_{\exists \forall}F$ they are free factors of $F$ by Theorem \ref{thm-elemantary-free-factor} and Theorem \ref{thm-free-factor}. Therefore $h_2 \circ h_1^{-1}$ can be extended to an automorphism of $F$ as required,  because $h_1(C)$ and $h_2(C)$ have the same rank. \qed
 
 \smallskip
 We conclude this section with the following proposition  of independent interest. 
 
 \begin{prop}Let $F$ be a nonabelian free group of finite rank  and let $\bar a$  be a tuple of $F$ such that $F$ is freely indecomposable relative to the subgroup generated by $\bar a$. Let $\bar s$ be a basis of $F$. Then there exists an universal formula $\varphi(\bar x)$ such that $F \models \varphi(\bar s)$ and such that for any endomorphism $f$ of $F$,  if $F \models \varphi(\bar s)$ and $f$ fixes $\bar a$ then $f$ is an automorphism.  In particular $(F, \bar a)$ is a prime model of the theory $Th(F, \bar a)$. 
 \end{prop}
 
 \proof  We claim that $\bar a$ is rigid. Suppose for a contradiction that $\bar a$ is not rigid and let $F=A*B$ be a nontrvial free decomposition such that $A$ contains a tuple $\bar c$ with $tp_{\exists \forall}^F(\bar a)=tp_{\exists \forall}^F(\bar c)$. By Theorem \ref{thm-elemantary-free-factor}, we conclude that $tp^F(\bar a)=tp^F(\bar c)$,  and by Proposition \ref{prop-princip},  there is an automorphism $\sigma$ sending $\bar c$ to $\bar a$. Hence $\sigma(A)$ is a free factor containing $\bar a$ and thus $F$ is freely indecomposable relative to the subgroup generated by $\bar a$; which is a contradiction to the hypothesis of the proposition.

 Hence,  by Proposition \ref{prop-key2},  there exists an universal  formula $\varphi(\bar y)$ such that $F \models \varphi(\bar s)$ and   such that for any $f\in Hom(F|\bar a,F|\bar a)$,  if $F \models \varphi(f(\bar s))$ then $f$ is an embedding and by Proposition \ref{prop1},  we conclude that $f$ is an automorphism. 
 
 Now the proof of the fact that $(F, \bar a)$ is a prime model of the theory $Th(F,\bar a)$ proceed in a similar way to that of Lemma \ref{main-lem3}. \qed
 
 \section{Dealing with parameters}

 In this section we show Theorem \ref{thm1} and \ref{thm2} for arbitrary  $P$. We reduce the problem to Proposition \ref{prop-dicho-exit} and Proposition \ref{prop-princip}  by using the definable closure and the existential definable closure. We recall the following definition. 
 
 \begin{defn} Let $G$ be a group and $P \subseteq G$. The \textit{definable closure} (resp. \textit{existential definable closure})  of $A$, denoted  $dcl(P)$ (resp. $dcl^{\exists}(P)$),  is the set of elements $g \in G$  such that there exists a formula  (resp. an existential formula)  $\phi(x)$ with parameters from $P$  such that $G \models \phi(g)$ and $g$ is the unique element satisfying $\phi$.  
 \end{defn}

We see that for any $P \subseteq G$, $dcl(P)$ and $dcl^{\exists}(P)$ are subgroups of $G$. In a furthcoming paper \cite{ould-dcl},  we answer a question of Z. Sela about the definable and the algebraic closure. We will use the following theorem of that paper. 

\begin{thm} \cite{ould-dcl} \label{thm-dcl}Let $F$ be a nonabelian free group of finite rank and let $P \subseteq F$. Then $dcl(P)$ and $dcl^\exists(P)$ are finitely generated and their rank is bounded by the rank of $F$.  \qed
\end{thm} 

Now we show the following simple lemma. 

\begin{lem} \label{lem-dcl}Let $G$ be a group and $P \subseteq G$. Let $\bar a$ $\bar b$ be  tuples from $G$.  

$(1)$ If  $tp^G(\bar a|P)=tp^G(\bar b|P)$ then $tp^G(\bar a|dcl(P))=tp^G(\bar b|dcl(P))$.

$(2)$ Similarly if $tp_{\exists}^G(\bar a|P)=tp_{\exists}^G(\bar b|P)$ then $tp_\exists ^G(\bar a|dcl^\exists(P))=tp_\exists ^G(\bar b|dcl^\exists(P))$.
\end{lem}

\proof $(1)$ Let $\psi(\bar x; y_1, \dots, y_n)$ be a formula such that $\psi(\bar x; d_1, \dots, d_n) \in tp^G(\bar a|dcl(P))$ where $d_i \in dcl(P)$ for all $i$. For every $i$, there exists a formula $\phi_i(y)$ with parameters from $P$ such that $d_i$ is the unique element satisfying $\phi_i$. Since 
$$
G \models \exists y_1, \dots, \exists y_n (\psi(\bar a; y_1, \dots, y_n) \wedge \bigwedge_{1 \leq i \leq n}\phi_i(y_i)),
$$
we find $g_1, \dots,g_n$ in $G$ such that 
$$
G \models  (\psi(\bar b; g_1, \dots, g_n) \wedge \bigwedge_{1 \leq i \leq n}\phi_i(g_i)),
$$
and thus we must have $g_i=d_i$ for all $i$. Therefore $\psi(\bar x; d_1, \dots, d_n) \in tp^G(\bar b|dcl(P))$.  Thus $tp^G(\bar a|dcl(P))\subseteq tp^G(\bar b|dcl(P))$ and by symmetry we conclude that $tp^G(\bar a|dcl(P))=tp^G(\bar b|dcl(P))$ as required.

$(2)$ The proof is similar to that of $(1)$ and it is left to the reader. \qed

\smallskip
\noindent
\textbf{Proof of Theorem \ref{thm1} } Let $F$ be a nonabelian group of finite rank. Let  $\bar a, \bar b \in F^n$ and let  $P \subseteq F$ such that $tp^F(\bar a|P)=tp^F(\bar b|P)$. 

By Lemma \ref{lem-dcl},  $tp^F(\bar a|dcl(P))=tp^F(\bar b|dcl(P))$,  and by Theorem \ref{thm-dcl} $dcl(P)$ is finitely generated. Let $\bar d$ be a basis of $dcl(P)$. We notice that $P \subseteq dcl(P)$. 

Then $tp^F(\bar a, \bar d )=tp^F(\bar b, \bar d )$ and thus there exists an automorphism $\sigma$ sending $\bar a$ to $\bar b$ and fixing $\bar d$ by Proposition \ref{prop-princip}. Thus in particular $\sigma$ fixes $P$. \qed

\smallskip
\noindent
\textbf{Proof of Theorem \ref{thm2}}. Similar to that of Theorem \ref{thm1} and the details are left to the reader. \qed

\section{Two-generated  torsion-free hyperbolic groups}

In this section we give the proof of Theorem \ref{thm3} and Corollary \ref{cor1}. We see that  Corollary \ref{cor1} is a mere consequence of Theorem \ref{thm3} and Lemma \ref{main-lem3}. It remains to show Theorem \ref{thm3}.  We first show the following lemma. 

\begin{lem} \label{lem-rigid-two} Let $\G$ be a torsion-free hyperbolic group and let $\G_0$ be a  subgroup of $\G$. Suppose that $\G=\<H, t| U^t=V\>$, where $U$ and $V$ are cyclic malnormal subgroups of $H$. If $\G_0$ is rigid then it is elliptic in the precedent splitting. 
\end{lem}

\proof  Suppose for a contradiction that $\G_0$ is not elliptic. Since $\G_0$ is freely indecomposable, $\G_0$ admits a cyclic splitting $\Lambda$ inherited by the given  splitting of $\G$. This splitting is non-principal because $\G_0$ is rigid. It follows that the graph corresponding to $\Lambda$ is a tree,  as otherwise $\G_0$ can be written as an HNN-extension,  contradicting again the rigidity of $\G_0$. In particular, $\G_0$ is an iterated amalgamated free product.

If each vertex group of $\Lambda$ is abelian, by the transivity of the commutation, $\G_0$ itself is abelian; a contradiction. 

Let $A_0$ be a nonabelian vertex group. We claim that each vertex group connected to $A_0$ is cylic. Let $V_0$ be the vertex group corresponding to $A_0$ and let $V_1$ be another vertex connected to $V_0$ by $e$. Let $\Lambda'$ be the graph obtained by deleting $e$ from $\Lambda$.   Then $\G_0=L_1*_{a=b}L_2$,  where $L_1$ and $L_2$ are the fundamental group of the connected components of $\Lambda'$. Then $A_0 \leq L_1$ or $A_0 \leq L_2$ and without loss of generality we assume that $A_0 \leq L_1$.  Hence $L_1$ is nonabelian. If $L_2$ is nonabelian, then $\G_0$ admits a principal cyclic splitting; a contradiction. Therefore, $L_2$ is abelian and thus cyclic.  The vertex group corresponding to $V_1$ is contained in $L_2$ and thus cyclic as claimed. 

Let $B_0$ be a vertex group corresponding to a vertex $V_1$ connected to $V_0$ by $e$. Since the splitting of $\G_0$ is inherited from that of $\G$,  the fundmental group of the graph of groups consisted of $V_0, V_1$ and $e$ is of the form $L=A^x*_{a^x=b^y}B^y$ where $A, B \leq H$ and $A_0=A^x, B_0=B^y$, $x,y \in \G$. 

We are going to show that $L$ is elliptic; that is $L$ is in a conjugate of $H$.  We have $a,b \in H$ and $a=(xy^{-1})b(yx^{-1})$ and $L^{x^{-1}}=A*_{a=b^{yx^{-1}}}B^{yx^{-1}}$. If $yx^{-1} \in H$ then $L$ is elliptic as claimed. So we suppose that $yx^{-1} \not \in H$.

Observe that $U \cap V^h =1$ for any $h \in H$; as otherwise $\G$ will contains $\mathbb Z^2$ which is a conradiction with its hyperbolicity. By Lemma \ref{lem1}, one of the following cases holds:

\smallskip
$(1)$   $b= {u^p}^{\gamma}$,     $a= {v^p}^{\delta}$,  $yx^{-1}=\gamma^{-1} t \delta$,  where $p \in \mathbb Z$ and $\gamma, \delta \in H$; 

\smallskip
$(2)$   $b= {v^p}^{\gamma}$,     $a= {u^p}^{\delta}$,  $yx^{-1}=\gamma^{-1} t^{-1} \delta$,  where $p \in \mathbb Z$ and $\gamma, \delta \in H$. 

\smallskip
where we have assumed that $U$ and $V$ are generated by $u$ and $v$ respectively. We treat only the case $(1)$,  the other case being similar.  Since $B$ is cyclic and $b= {u^p}^{\gamma}$ and $U$ is malnormal,  we get $B=\<{u^q}^{\gamma}\>$ for some $q \in \mathbb Z$. Therefore $$B=\delta^{-1}.t^{-1}. \gamma \<\gamma^{-1}u^q\gamma\>\gamma^{-1}.t.\delta=\<v^q\>^\delta \leq H,$$
and thus $L$ is elliptic as claimed. 

Let $\Lambda''$ be the graph of groups obtained by collapsing $e$. Then $\Lambda''$ has less vertices than $\Lambda$. Proceeding by induction on the number of vertices, we conclude that $\G_0$ is elliptic; a final contradiction. \qed

\smallskip
We are now in a position to prove Theorem \ref{thm3}. Let $\G_1 \leq \G_2 \leq \dots \leq \G_n=\G$ be a sequence given by Theorem \ref{thm-structure-two}, where $\G_1$ is rigid. Since $\G_1$ is rigid, by Lemma \ref{main-prop},  there exists a finite subset $S \subseteq \G_1 \setminus \{1\}$ such that for any endomorphism $\varphi$  of $\G$ if $1 \not \in \varphi(S)$ then $\varphi$ is one-to-one in restriction to $\G_n$.

 Let $\varphi $ be an endomorphism of $\G$ such that $1 \not \in \varphi(S)$ and let us show that $\varphi$ is an automorphism.

We have $\varphi(\G_1) \leq \G$ and $\G=\<\G_{n-1}, t| A^t=B\>$.  Now $\varphi(\G_1)$ is isomorphic to $\G_1$ and thus rigid. By Lemma \ref{lem-rigid-two},     $\varphi(\G_1)$ is elliptic in the above splitting; that is $\varphi(\G_1)$ is in a conjugate of $\G_{n-1}$.  Using a similar argument and proceeding by induction,  $\varphi(\G_1)$ is in a conjugate of $\G_1$. 

Let $g \in \G$ such that $\varphi(\G_1) \leq g\G_1g^{-1}$ and let $\tau_g(x)=x^g$. Therefore we have $\tau_g \circ \varphi(\G_1) \leq \G_1$.  Since $\G_1$ is co-hopfian and $\tau_g \circ \varphi$ is one-to-one in restriction to $\G_1$, we conclude that $\tau_g \circ \varphi(\G_1)=\G_1$.

Set $\phi=\tau_g \circ \varphi$. We show by induction on $i$ that the restriction of $\phi$ to $\G_{i}$ is an automorphism of $\G_i$.  Write 
$$
\G_{i+1}=\<\G_{i}, t_{i}| A_{i}^{t_{i}}=B_{i}\>, \; A_i=\<a_i\>, \; B_i=\<b_i\>. 
$$

We claim that $\phi(t_1) \in \G_{2}$. If $n=2$  clearly $\phi(t_1) \in \G_{2}$. Hence we suppose  that $n \geq 3$. 

Let us first prove that $\phi(t_1) \in \G_{n-1}$. Suppose for  a contradiction that $\phi(t_1) \not \in \G_{n-1}$. We have  $$\phi(t_{1})^{-1}\phi(a_1)\phi(t_1)=\phi (b_1)\hbox{ and }\phi(a_1), \phi(b_1) \in \G_1 \leq \G_{n-1}. $$

Observe that $U \cap V^h =1$ for any $h \in H$; as otherwise $\G$ will contains $\mathbb Z^2$ which is a conradiction with its hyperbolicity. According to Lemma \ref{lem1}, one of the following cases holds: 

\smallskip
$(1)$   $\phi(a_1)=\gamma^{-1} a_{n-1}^{p}\gamma$,     $\phi(b_1)= \delta^{-1}{b_{n-1}^p}\delta$,  $\phi(t_1)=\gamma^{-1} t_{n-1} \delta$,  where $ p \in \mathbb Z$ and $\gamma, \delta \in \G_{n-1}$. 

\smallskip
$(2)$   $\phi(a_1)=\gamma^{-1} {b_{n-1}^p} \gamma$,     $\phi(b_1)=\delta^{-1} {a_{n-1}^p}\delta$,  $\phi(t_1)=\gamma^{-1} t_{n-1}^{-1} \delta$,  where $p \in \mathbb Z$ and $\gamma, \delta \in \G_{n-1}$. 

\smallskip
Let us treat the case $(1)$,  the  case $(2)$  can be treated similarly.   We first show that $p=\pm 1$. We have $a_{n-1}^\gamma \in C_{\G_{n-1}}(\phi(a_1))$.  According to  \cite[Theorem 3.2(i)]{ould-e.c.groups}, $C_{\G_{n-1}}(\phi(a_1))=C_{\G_{n-2}}(\phi(a_1))$.  A repeated application of \cite[Theorem 3.2(i)]{ould-e.c.groups}, gives $C_{\G_{n-1}}(\phi(a_1))=C_{\G_1}(\phi(a_1))$.  Therefore $a_{n-1}^\gamma \in C_{\G_1}(\phi(a_1))$.  Since the restriction of $\phi$ to $\G_1$ is an automorphism,  we find $c \in \G_1$ such that $\phi(c)=a_{n-1}^\gamma$ and $a_1=c^p$.  Since $a_1$ is root-free, we conclude finally that $p=\pm 1$ as claimed.

We rewrite now $\G$ as follows 
 $$
 \G=\<\G_{n-1}, s| s^{-1} \phi(a_1)s=\phi(b_1)\>,
 $$
 where $s=\phi(t_1)$.  We also have 
 $$
 \G_{n-1}=\<\G_{n-2}, t_{n-2}| A_{n-2}^{t_{n-2}}=B_{n-2}\>, \; \phi(a_1), \phi(b_1)  \in \G_{n-2}. 
 $$
 
Hence $\G$ admits a principal cyclic splitting with more than one edge; a contradiction with \cite[Theorem A]{kapovich-two}.  Therefore $\phi(t_1) \in \G_{n-1}$ as claimed.  
 
 Using a similar argument and proceeding  by induction we conclude that  $\phi(t_1) \in \G_{2}$. In particular $\phi(\G_{2}) \leq \G_{2}$. 
 
 Clearly $ \phi(t_1) \not \in \G_1$;  otherwise $a_1$ and $b_1$ are conjugate in $\G_1$ and thus $\G_{2}$ contains $\mathbb Z^2$ contradicting its hyperbolicity. Hence, by Lemma \ref{lem1}, one of the following cases holds 
 
\smallskip
$(1)$   $\phi(a_1)=\gamma^{-1} a_{1}^{p}\gamma$,     $\phi(b_1)= \delta^{-1}{b_{1}^p}\delta$,  $\phi(t_1)=\gamma^{-1} t_{1} \delta$,  where $ p \in \mathbb Z$ and $\gamma, \delta \in \G_{1}$. 

\smallskip
$(2)$   $\phi(a_1)=\gamma^{-1} {b_{1}^p} \gamma$,     $\phi(b_1)=\delta^{-1} {a_{1}^p}\delta$,  $\phi(t_1)=\gamma^{-1} t_{1}^{-1} \delta$,  where $p \in \mathbb Z$ and $\gamma, \delta \in \G_{1}$. 

\smallskip

Let us treat the case $(1)$,  the  case $(2)$  being similar. Proceeding as above, we have $p=\pm 1$.  Again as before,  we rewrite $\G_2$ as
$$
\G_2=\<\G_{1}, s| s^{-1} \phi(a_1)s=\phi(b_1)\>,
$$
where $s=\phi(t_1)$.  Hence,  we get $\phi(\G_{2})=\G_{2}$ and in particular the restriction of $\phi$ to $\G_2$  is an automorphism of $\G_2$.

Applying the same argument and proceeding by induction,  we conclude that for every $i$ the restriction of $\phi$ to $\G_i$ is an automorphism of $\G_i$.  In  particular $\phi$ is an automorphism of $\G$ as well as  $\varphi$. \qed

\section{Remarks}

$(1)$ We note   that a non-free two-generated torsion-free hyperbolic group is not necessarily rigid. Here an example.  Let $F=\<a,b|\>$ be the free group of rank $2$ and let $r \in F$ satisfying the following properties:

\smallskip
$(i)$ $r$ is root-free, is cyclically reduced  and its lenght is greater than $6$; 

$(ii)$   the symmetrized set generated by $r$ satisfies $C'(1/8)$. 

\smallskip
Let $\G=\<a,b|r=1\>$. By \cite[Theorem 5.4, V]{LyndonSchupp77}, $a$ and $b$ are not conjugate in $\G$. It follows in particular that any power of $a$ is not conjugate to any power of $b$.  We see also that $\<a|\>$ and  $\<b|\>$  are  malnormal in $\G$. Hence the HNN-extension $L=\<\G, t | a^t=b\>$ is conjugately seperated in the sense of \cite{Kharlam-free-construction}. Since $\G$ is torsion-free  and hyperbolic, by \cite[Corollary 1]{Kharlam-free-construction}, $L$ is a torsion-free hyperbolic group.   Hence $L$ is a non-free two-generated torsion-free hyperbolic group which admits an essential cyclic splitting; and thus  $L$ is not rigid.

$(2)$ It is noted in \cite{Pillay-fork} that  a nonabelian free group is connected. Hence,  one may ask if this is still true for nonabelian torsion-free hyperbolic groups. Recall that a group $G$ is said to be \emph{connected},  if $G$ is without definable subgroup of finite index.

\begin{prop} A noncyclic torsion-free hyperbolic group is connected. 
\end{prop}

Recall that a definable subset $X$ of $G$ is said to be \emph{right generic},  if there exist $g_1, \cdots, g_n \in G$ such that $G= g_1 X \cup \cdots \cup g_nX$. Left generic definable subsets are defined anagousely. Now we show the following lemma. 

\begin{lem} \label{lem-gen} Let $G$ be a group and suppose that $G$ satisfies: if $X$ and $Y$ are right generic sets then $X \cap Y \neq \emptyset$. Then $G$ is connected. 
\end{lem}

\proof

If $H$ is a definable subgroup of finite index,  then $G= g_1 H \cup \cdots \cup g_nH$. Then any $g_iH$ is  right generic because
$$
G=(g_1g_i^{-1})g_iH \cup \cdots (g_ng_i^{-1})g_iH,
$$
and therefore for any $i,j$,  $g_iH \cap g_jH \neq \emptyset$ and thus we must have $G=H$. \qed

\smallskip
For a group $G$, we denote by $G[a]$ the group $G*\mathbb Z$ where $a$ is a generating element of $\mathbb Z$.  The following lemma is a slight rafinement of an observation of B .Poizat. 

\begin{lem} \label{lem-con}Let $G$ be a group and suppose that  $G \preceq G[a]$. If $\phi(G)$ is a  right generic subset of $G$,  then $a
\in \phi(G[a])$. In particular $G$ is connected. \end{lem}

\proof

Suppose that $G =g_1  \phi(G) \cup \cdots \cup g_n \phi(G)$.
Since $G \preceq G[a]$ we get also $G[a]=g_1 \phi(G[a]) \cup
\cdots \cup g_n \phi(G[a])$. Therefore for some $i$, $g_i^{-1}a
\in \phi(G[a])$. Since there exists an automorphism of $G[a]$, fixing pointwise $G$, 
which sends $g_i^{-1}a$ to $a$,  we get $a \in \phi(G[a])$.

Since $G \preceq G[a]$,  it follows that if $X$ and $Y$ are right generic subsets,  then $X \cap Y \neq \emptyset$.  Therefore by   Lemma \ref{lem-gen}, $G$ is connected. \qed

\smallskip
Let $\G$ be a torsion-free hyperbolic group. In \cite{sela-hyp},  Z. Sela shows that if $\G$ is not elementary equivalent to a free group, then $\G$ has a smallest elementary subgroup, denoted by $EC(\G)$, called the \textit{elementary core} of $\G$.  For the complete definition, we refer the reader to \cite[Definition 7.5]{sela-hyp}. We need the following properties which follows from \cite{sela-hyp}. 

\begin{fact} The elementary core satisfies the following properties: 

$(1)$ $EC(\G)=1$ if and only if $\G$ is elementary equivalent to a nonabelian free group. 

$(2)$ $EC(\G)=EC(\G*\mathbb Z)$ and $EC(EC(\G))=EC(\G)$. 

$(3)$ If $EC(\G)\neq 1$ then $EC(\G) \preceq \G$. 
\end{fact}

By $(2)$ and $(3)$ we get $EC(\G) \preceq EC(\G)*\mathbb Z$. We conclude, by Lemma \ref{lem-con},  that $EC(\G)$ is connected and thus by elementary equivalence,  $\G$ is connected.

\bigskip
\noindent Abderezak OULD HOUCINE, \\ 
Universit\'e  de Mons, Institut de Mathmatique, B\^atiment Le Pentagone, 
Avenue du Champ de Mars 6, B-7000 Mons, Belgique.\\ 

\noindent Universit\'e  de Lyon; Universit\'e Lyon 1; INSA de Lyon, F-69621; Ecole Centrale 
de Lyon; CNRS, UMR5208, Institut Camille Jordan, 43 blvd du 11 novembre 
1918, F-69622 Villeurbanne-Cedex, France.  \\

\noindent \textit{E-mail}:\textrm{ould@math.univ-lyon1.fr}

\bibliographystyle{alpha}
\bibliography{biblio}

\def\cprime{$'$} \def\cprime{$'$} \def\cprime{$'$}
\begin{thebibliography}{BMR99}

\bibitem[Bel07]{Beleg}
I.~Belegradek.
\newblock Aspherical manifolds with relatively hyperbolic fundamental groups.
\newblock {\em Geom. Dedicata}, 129:119--144, 2007.

\bibitem[BMR99]{Baum-Rem-algebraic}
G.~Baumslag, A.~Myasnikov, and V.~Remeslennikov.
\newblock Algebraic geometry over groups. {I}. {A}lgebraic sets and ideal
  theory.
\newblock {\em J. Algebra}, 219(1):16--79, 1999.

\bibitem[CG05]{Champ-Guirardel}
C.~Champetier and V.~Guirardel.
\newblock Limit groups as limits of free groups.
\newblock {\em Israel J. Math.}, 146:1--75, 2005.

\bibitem[Del91]{delzant-two}
T.~Delzant.
\newblock Sous-groupes \`a deux g\'en\'erateurs des groupes hyperboliques.
\newblock In {\em Group theory from a geometrical viewpoint ({T}rieste, 1990)},
  pages 177--189. World Sci. Publ., River Edge, NJ, 1991.

\bibitem[GW07]{groves-2007}
D.~Groves and H.~Wilton.
\newblock Conjugacy classes of solutions to equations and inequations over
  hyperbolic groups, 2007.

\bibitem[Hod93]{Hodges(book)93}
W.~Hodges.
\newblock {\em Model theory}, volume~42 of {\em Encyclopedia of Mathematics and
  its Applications}.
\newblock Cambridge University Press, Cambridge, 1993.

\bibitem[JOH04]{ould-e.c.groups}
E.~Jaligot and A.~Ould~Houcine.
\newblock Existentially closed {CSA}-groups.
\newblock {\em J. Algebra}, 280(2):772--796, 2004.

\bibitem[KM98]{Kharlam-free-construction}
O.~Kharlampovich and A.~Myasnikov.
\newblock Hyperbolic groups and free constructions.
\newblock {\em Trans. Amer. Math. Soc.}, 350(2):571--613, 1998.

\bibitem[KM06]{Kharlam-Mya-free-gps}
Olga Kharlampovich and Alexei Myasnikov.
\newblock Elementary theory of free non-abelian groups.
\newblock {\em J. Algebra}, 302(2):451--552, 2006.

\bibitem[KW99]{kapovich-two}
I.~Kapovich and R.~Weidmann.
\newblock On the structure of two-generated hyperbolic groups.
\newblock {\em Math. Z.}, 231(4):783--801, 1999.

\bibitem[LS77]{LyndonSchupp77}
R.~C. Lyndon and P.~E. Schupp.
\newblock {\em Combinatorial group theory}.
\newblock Springer-Verlag, Berlin, 1977.
\newblock Ergebnisse der Mathematik und ihrer Grenzgebiete, Band 89.

\bibitem[Mar02]{Marker}
D.~Marker.
\newblock {\em Model Theory : An Introduction}.
\newblock Springer-Verlag, New york, 2002.
\newblock Graduate Texts in Mathematics.

\bibitem[Nie03a]{nies-free}
A.~Nies.
\newblock Aspects of free groups.
\newblock {\em J. Algebra}, 263(1):119--125, 2003.

\bibitem[Nie03b]{Nies-QFA}
A.~Nies.
\newblock Separating classes of groups by first-order sentences.
\newblock {\em Internat. J. Algebra Comput.}, 13(3):287--302, 2003.

\bibitem[Nie07]{Nies-QFA2}
A.~Nies.
\newblock Comparing quasi-finitely axiomatizable and prime groups.
\newblock {\em J. Group Theory}, 10(3):347--361, 2007.

\bibitem[OH]{ould-csa-subgroup}
A.~Ould~Houcine.
\newblock Subgroup theorem for valuated groups and the csa property.
\newblock {\em To appear in Journal of Algebra}.

\bibitem[OH07]{Ould-equa}
A.~Ould~{H}oucine.
\newblock Limit {G}roups of {E}quationally {N}oetherian {G}roups.
\newblock In {\em Geometric group theory}, Trends Math., pages 103--119.
  Birkh\"auser, Basel, 2007.

\bibitem[OH08]{Ould-csa-super}
A.~Ould~{H}oucine.
\newblock On superstable {CSA}-groups.
\newblock {\em Annals of Pure and Applied Logic}, 154(1):1--7, 2008.

\bibitem[OH09a]{ould-cb-limit}
A.~Ould~{H}oucine.
\newblock Note on the {C}antor-{B}endixson rank of limit groups.
\newblock {\em To appear in Communications in Algebra}, 2009.

\bibitem[OH09b]{ould-rank}
A.~Ould~{H}oucine.
\newblock On finitely generated models of theories with few finitely generated
  models.
\newblock {\em Submitted}, 2009.

\bibitem[OHV10]{ould-dcl}
A.~Ould~{H}oucine and D.~Vallino.
\newblock Algebraic and definable closure in free groups.
\newblock {\em In preparation}, 2010.

\bibitem[Ol{\cprime}80]{Ol-monster}
A.~Yu. Ol{\cprime}shanski{\u\i}.
\newblock An infinite group with subgroups of prime orders.
\newblock {\em Izv. Akad. Nauk SSSR Ser. Mat.}, 44(2):309--321, 479, 1980.

\bibitem[Ol{\cprime}91]{book-olshanski}
A.~Yu. Ol{\cprime}shanski{\u\i}.
\newblock {\em Geometry of defining relations in groups}, volume~70 of {\em
  Mathematics and its Applications (Soviet Series)}.
\newblock Kluwer Academic Publishers Group, Dordrecht, 1991.
\newblock Translated from the 1989 Russian original by Yu. A. Bakhturin.

\bibitem[OS06]{Oger-Sabbagh}
F.~Oger and G.~Sabbagh.
\newblock Quasi-finitely axiomatizable nilpotent groups.
\newblock {\em J. Group Theory}, 9(1):95--106, 2006.

\bibitem[OT00]{john-turner}
J.~C. O'Neill and E.~C. Turner.
\newblock Test elements and the retract theorem in hyperbolic groups.
\newblock {\em New York J. Math.}, 6:107--117 (electronic), 2000.

\bibitem[Per08]{chloe-these}
C.~Perin.
\newblock Plongements \'el\'ementaires dans un groupe hyperbolique sans
  torsion.
\newblock {\em Th\`ese de doctorat, Universit\'e de Caen/Basse-Normandie},
  2008.

\bibitem[Pil08]{Pillay-fork}
Anand Pillay.
\newblock Forking in the free group.
\newblock {\em J. Inst. Math. Jussieu}, 7(2):375--389, 2008.

\bibitem[Pil09]{pillay-primitive}
Anand Pillay.
\newblock On genericity and weight in the free group.
\newblock {\em Proc. Amer. Math. Soc.}, 137(11):3911--3917, 2009.

\bibitem[PS10]{perin-hom}
C.~Perin and R.~Sklinos.
\newblock Homogeneity in the free group.
\newblock {\em Preprint}, 2010.

\bibitem[RS94]{rips-sela-structure-rigidity}
E.~Rips and Z.~Sela.
\newblock Structure and rigidity in hyperbolic groups. {I}.
\newblock {\em Geom. Funct. Anal.}, 4(3):337--371, 1994.

\bibitem[Sel97]{Sela-hopf}
Z.~Sela.
\newblock Structure and rigidity in ({G}romov) hyperbolic groups and discrete
  groups in rank {$1$} {L}ie groups. {II}.
\newblock {\em Geom. Funct. Anal.}, 7(3):561--593, 1997.

\bibitem[Sel05]{Sela-Dio-V1}
Z.~Sela.
\newblock Diophantine geometry over groups. {$\rm V\sb 1$}. {Q}uantifier
  elimination. {I}.
\newblock {\em Israel J. Math.}, 150:1--197, 2005.

\bibitem[Sel06a]{Sela-Dio-V2}
Z.~Sela.
\newblock Diophantine geometry over groups. {${\rm V}\sb 2$}. {Q}uantifier
  elimination. {II}.
\newblock {\em Geom. Funct. Anal.}, 16(3):537--706, 2006.

\bibitem[Sel06b]{Sela-Dio-VI}
Z.~Sela.
\newblock Diophantine geometry over groups. {VI}. {T}he elementary theory of a
  free group.
\newblock {\em Geom. Funct. Anal.}, 16(3):707--730, 2006.

\bibitem[Sel09]{sela-hyp}
Z.~Sela.
\newblock Diophantine geometry over groups. {VII}. {T}he elementary theory of a
  hyperbolic group.
\newblock {\em Proc. Lond. Math. Soc. (3)}, 99(1):217--273, 2009.

\end{thebibliography}
\end{document}